%% file: swept2d-paper.tex
\documentclass[review]{elsarticle}
\usepackage{amsmath, amssymb}
\usepackage{graphicx}
\usepackage{caption}
\usepackage{subcaption}
\usepackage{hyperref}
\usepackage{lineno,hyperref}
\usepackage{listings}
\usepackage{xcolor}
\usepackage[ruled]{algorithm2e}
\usepackage{txfonts}
\usepackage{verbatim}
\modulolinenumbers[5]

\journal{Journal of Computational Physics}

\begin{document}

\begin{frontmatter}

\title{The swept rule for breaking the latency barrier in time advancing two-dimensional PDEs}


\author[address1]{Maitham Alhubail\corref{mycorrespondingauthor}}
\cortext[mycorrespondingauthor]{Corresponding author}
\ead{hubailmm@mit.edu}

\author[address1]{Qiqi Wang}
\ead{qiqi@mit.edu}
\ead[url]{www.engineer-chaos.blogspot.com}

\author[address2]{John Williams}
\ead{jrw@mit.edu}

\address[address1]{Aerospace Computational Design Laboratory, MIT, 77 Massachusetts Ave., Cambridge, MA 02139, USA}

\address[address2]{MIT Geonumerics Group, MIT, 77 Massachusetts Ave., Cambridge, MA  02139, USA}

\begin{abstract}
This article describes a method to accelerate parallel, explicit time integration of two-dimensional unsteady PDEs. The method is motivated by our observation that network latency, not bandwidth or computing power, often limits how fast PDEs can be solved in parallel. The method is called the swept rule of space-time domain decomposition. Compared to conventional, space-only domain decomposition, it communicates similar amount of data, but in fewer messages. The swept rule achieves this by decomposing space and time among computing nodes in ways that exploit the domains of influence and the domain of dependency, making it possible to communicate once per many time steps with no redundant computation.  By communicating less often, the swept rule effectively breaks the latency barrier, advancing on average more than one time step per ping-pong latency of the network. The article describes the algorithms, presents simple theoretical analysis to the performance of the swept rule in two spatial dimensions, and supports the analysis with numerical experiments.  
\end{abstract}


\end{frontmatter}


\include{introduction}
\include{technical}

\include{analysis}

\include{interface}
\include{experiments}
\include{conclusion}

\section*{References}

\bibliography{mybibfile}

\end{document}

%% file: introduction.tex
\section{Introduction}

There is a strong demand to numerically solve partial differential
equations (PDEs) faster by using more parallel processors.  Today, this
demand is unmet for small and long unsteady simulations.
A simulation is small if the PDE is discretized spatially
on a few million grid points.  In three spatial dimensions, this
means that each dimension may be discretized into just a few hundreds
grid points.  
A simulation is long if the PDE need to be time integrated for millions
of time steps.

Small and long simulations
arrise, for example, in aero-thermal design of turbomachinary components.
The simulation can be small when it is used to design a local feature of a
component, e.g., the trailing edge of a turbine blade.
It can nevertheless be long, because the slow thermodynamics
requires a long simulation time span, and the fast, unsteady aerodynamics
forces this long time span to be divided into millions of tiny time steps.
Other simulations can be small and long if it resolves tightly coupled
physical processes with similar spatial scales but disparate time scales.
Such simulations are difficult to parallelize, and even more difficult to
scale to many parallel processors.

Scaling of a parallel PDE solver is always limited.  When the solver runs
on a few processors, doubling the processors can half the run time.
As the number of processors further doubles, however, the percentage
reduction in run time starts to diminish.
As the number of processors reaches the scaling limit,
adding more processors no longer reduces the run time.

For a well-engineered parallel PDE solver, what limits the scaling is the
communication between computing nodes.  As the same computational task
is subdivided into more computing nodes, the nodes communicate more frequently
to exchange smaller batches of data.  Exchanging smaller batches of
data, however, is not always faster.  Even the smallest batch of data
takes a finite, predictable amount of time to be exchanged.  That amount of
time is the network latency.  It is a fundamental cause to the scaling limit.

This paper focuses on extending the scaling limit by circumventing
the network latency.  It achieves this through exchanging fewer, larger batches
of data between computing nodes.  Because of the network latency, how long it
takes to exchange a hundred or fewer bytes of data 
does not depend on the amount of data.  On most networks, exchanging
these bytes in a one shot is almost exactly twice as fast as
exchanging them in two batches, if the first batch must finish before the
second starts.  Latency affect the speed of data
exchange in batches up to about a megabyte.  For less than a megabyte of
data, exchanging it in a single batch is noticeably faster on most networks
than exchanging it into multiple batches.  Therefore, by communicating
similar amount of data in fewer batches, the algorithm described in this paper
delays the effect of network latency, thereby extending the scaling limit
of parallel PDE solvers.

The algorithm described in this paper, namely the swept rule, communicates
less often and incurs less network latency.  It achieves so by decomposing
space and time among computing nodes with respect to domains of influence
and domains of dependency.  In this aspect, it is inspired by the
Communication Avoiding (CA) algorithm \cite{ca1, ca2} by Demmel et al.
Unlike the CA
algorithm, however, the swept rule does not incur redundant computation.
The swept rule shares some similarity with diamond tiling for cache
optimization\cite{tilt1,tilt2,tilt3}.

The swept rule joins the ranks of space-time parallel methods
for solving PDEs \cite{50years},
other members of which include Parallel Full Approximation Scheme
in Space and Time (PFASST) \cite{pfasst}, Parareal \cite{parareal},
Space-Time Parallel Multigrid algorithms
\cite{mgrit}, Parallel Implicit Time-integration Algorithms (PITA),
Adjoint-based time-parallel integration \cite{adjoint}, and Parallel-in-time
algorithm for chaotic systems \cite{wang2013towards}.
The algorithm described in this paper can either be used alone, or in
conjunction with other space-time parallel methods to further delay the
scaling limit.

We previously introduced the swept rule in the context of unsteady PDEs in
a single spatial dimension\cite{swept}.  This article extends the swept rule
to PDEs in two spacial dimensions.  We do not assume that the reader has
read the one-dimensional swept rule paper; however, some fundamental
concepts are explained in that paper from a different angle or in more detail.
In this paper, we show how the swept rule works in two spatial dimensions,
and how it breaks the latency barrier via communicating less often.
The basic idea behind Swept rule in 2D is the same to that in 1D.
It decomposes space and time among computing nodes in ways that exploit
the domains of influence and the domain of dependency, making it possible
to communicate once per many time steps. The resulting algorithm enables
simulations to be solved significantly faster than what is possible with
spatial domain decomposition schemes typically found in today's PDE solvers.

The remainder of the paper is organized as follows.
Section 2 introduces the setting of the scheme developed in this paper, and under this setting, the discrete space-time of a PDE solver. Sections 3 and 4 introduce the Swept 2D components and how they work with each other to break the latency barrier.  Sections 5 and 6 analyzes the performance of Swept 2D and present our interface to utilizing our implementation of Swept 2D.  Sections 7 and 8 present the results of two experiments done using our implementation of Swept 2D with two PDE applications and finally, we conclude with a summary.

%% file: technical.tex
\section{Space-time Decomposition of a PDE Solver}

The swept rule described in this paper operates under the following set
of assumptions:
\begin{enumerate}
\item The spatial domain is discretized into a doubly-periodic, logically
      Cartesian grid.  Each point in the Cartesian grid has a pair of indices
      $(i,j)$, and is considered to have 8 neighbors, with indices
      $(i,j\pm1), (i\pm1,j), (i\pm1,j\pm1)$.
\item The scheme for advancing each time step
      computes a few quantities on every grid point at the next time step,
      using the quantities on its neighbors at the current time step.
      Alternatively, a stencil operation that needs access beyond its
      nearest neighbors should be decomposed into sub-steps that accesses
      only the nearest neighbor.
\end{enumerate}

Such decomposition is feasible for complex explicit schemes, and can
be automated\cite{stencil}.  It eliminates the need of accessing the neighbors
of neighbors.  Here we illustrate such decomposition with an example.
Consider the following scheme which requires two levels of neighbors.
\begin{equation} \label{complex_step}
    u^{n+1}_{i,j} = u^n_{i-2,j} + u^n_{i+2,j} + u^n_{i,j-2} + u^n_{i,j+2}
- 4 u^n_{i,j} \end{equation}
It can be decomposed into 2 sub-steps.
The first sub-step simply ``pushes'' its neighbor values to the second
sub-step.
\begin{equation} \label{substep1}
    u^{n+\frac{1}{2}}_{i,j} = \begin{bmatrix}
   u^{n}_{i,j} \\ u^{n}_{i-1,j}\\ u^{n}_{i+1,j}\\ u^{n}_{i,j-1}\\ u^{n}_{i,j+1}
\end{bmatrix}
\end{equation}
The second sub-step has access to $u^{n+\frac{1}{2}}$ on the neighbors, which
contains $u^{n}$ on the neighbors of neighbors, which were ``pushed'' by
the first sub-step.  Using $u^{n+\frac12}_{i,j;k}$ to denote the $k$th
component ($k=1,\ldots,5$) of the vector $u^{n+\frac12}_{i,j}$, the
second sub-step should be
\begin{equation} \label{substep2}
    u^{n+1}_{i,j} = u^{n+\frac12}_{i-1,j;2} + u^{n+\frac12}_{i+1,j;3} + u^{n+\frac12}_{i,j-1;4} + u^{n+\frac12}_{i,j+1;5}
                 - 4 u^{n+\frac12}_{i,j;1}
\end{equation}
These 2 sub-steps (\ref{substep1}-\ref{substep2}) computes the same quantities as (\ref{complex_step}).
But each sub-step only uses the immediate neighbors.  Eliminating the need
for further grid points can be achieved by decomposing into more sub-steps
\cite{stencil}.

\section{The Components of Swept Rule in 2D}

We break the Swept Rule in two dimensional domains (Swept 2D) into
four simple components.  In this section, we first describes how these
components work with each other to break the latency barrier.  We then
describe in detail how each of these components can be built.
These four components are upward pyramid, longitudinal bridge,
latitudinal bridge, and downward pyramid.
All these components are a subdomain in space-time,
whose boundaries follow the discrete domain of dependence and the domains
of influence.

\begin{figure}
	\centering
	\begin{subfigure}[b]{0.45\textwidth}
	\includegraphics[width=1.0\textwidth]{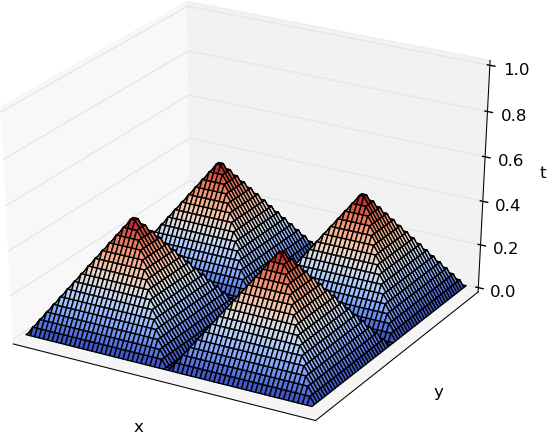}
	\caption{Starting with square partitions and performing parallel computation without any communication to form pyramids}
	\end{subfigure}
	\quad
	\begin{subfigure}[b]{0.45\textwidth}
	\includegraphics[width=1.0\textwidth]{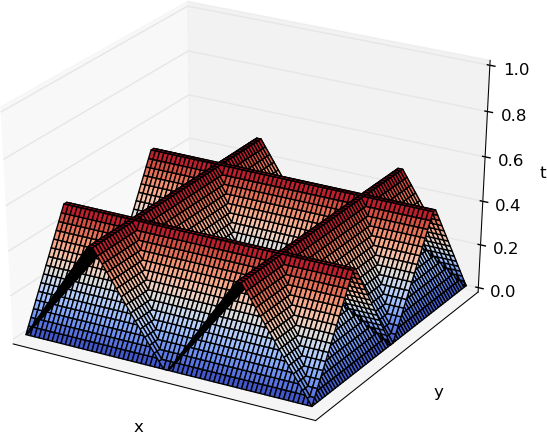}
	\caption{Building the bridges that connect the tips of the pyramids and filling the gaps between the built bridges}
	\end{subfigure}

	\begin{subfigure}[b]{0.45\textwidth}
	\includegraphics[width=\textwidth]{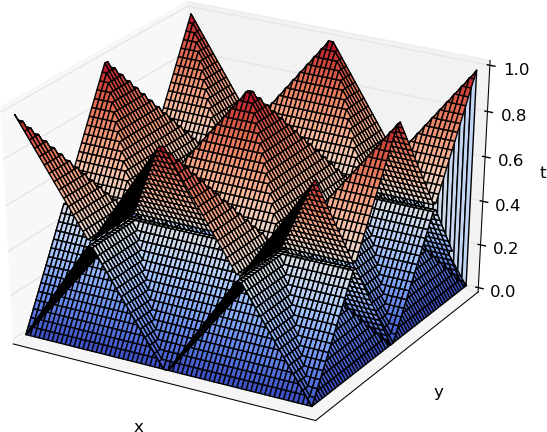}
	\caption{Growing the next set of pyramids\newline\newline}
	\end{subfigure}
	\quad
	\begin{subfigure}[b]{0.45\textwidth}
	\includegraphics[width=\textwidth]{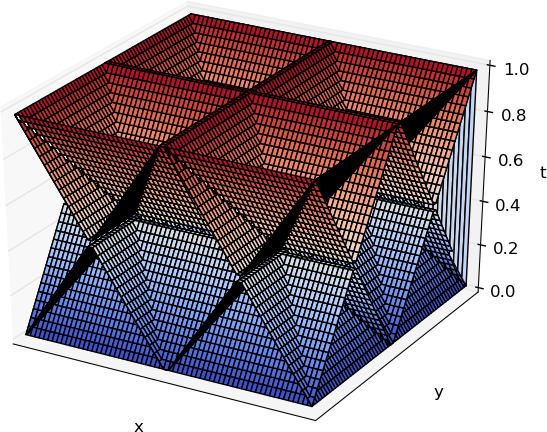}
	\caption{Completing the Swept 2D cycle by building the next set of bridges and filling the gaps between them}
	\end{subfigure}	
		
    \caption{The big picture of the Swept Rule in 2D.  
Because of periodic boundary condition, the half pyramids and quarter pyramids parts of full pyramids.}
 \label{bigpicture}
\end{figure}

Figure \ref{bigpicture} illustrates the idea of the Swept Rule in 2D.
To start, we partition the spatial domain into squares subdomains.
According to this partitioning, the initial condition of the PDE is
distributed among the processors.  Using the initial condition on a
subdomain, each processor builds an upward pyramid, which is
the first component of Swept 2D.  These upward pyramids are illustrated
in Fig \ref{bigpicture}(a).  After that, each processor sends data to
two neighboring processors, and receives data from two other neighboring
processors.  Once the data is received, each processor builds
one longitudinal bridge and one latitudinal bridge, and then sends and receives
data from its neighboring processors again.  The bridges, which are
the second and third components in Swept 2D, are illustrated in
Fig \ref{bigpicture}(b).  Finally, the processor fills
the gaps between the generated bridges with a downward pyramid, the fourth
and last component.  This completes what we call a half Swept 2D cycle.
Repeating this process, illustrated in Fig \ref{bigpicture}(c-d),
completes a full Swept 2D cycle.  The rest of this section explains the
details of each of the four components.

\subsection{Upward Pyramid}

The upward pyramid is a square-pyramid-shaped subdomain in the three-dimensional, discrete
space-time. The two spatial dimensions are discretized with a grid indexed by $(i,j)$; the time dimension
is discretized with time steps indexed by $k$.  Without loss of generality, denote the first time step in the upward pyramid as $0$,
last time step as $n/2-1$.  The base of the pyramid 
is a square subdomain of $n$ by $n$ grid points at time step $0$.  As the time step increases,
the cross section of the upward pyramid maintains a square shape, whose side length decreases
by 2 for every time step.  At the last time step, the cross section is a 2 by 2 square.

\begin{figure}
	\centering
	\begin{subfigure}[b]{0.3\textwidth}
	\includegraphics[width=1.0\textwidth]{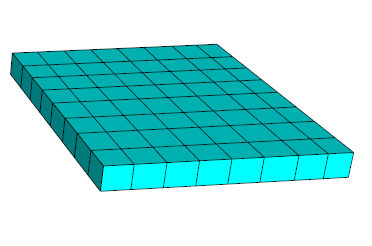}
	\caption{Level 0 (8x8)}
	\end{subfigure}
	\begin{subfigure}[b]{0.3\textwidth}
	\includegraphics[width=1.0\textwidth]{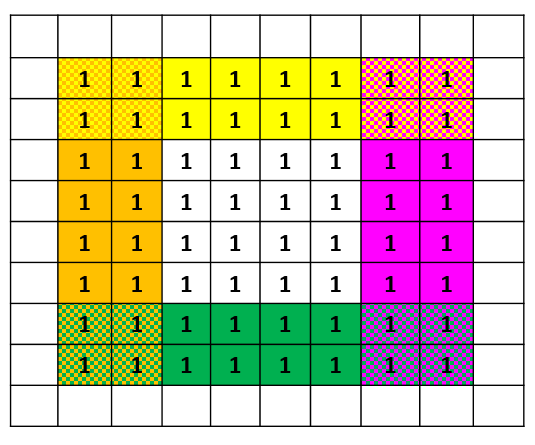}
	\caption{Populating the sides}
	\end{subfigure}
	\begin{subfigure}[b]{0.3\textwidth}
	\includegraphics[width=1.0\textwidth]{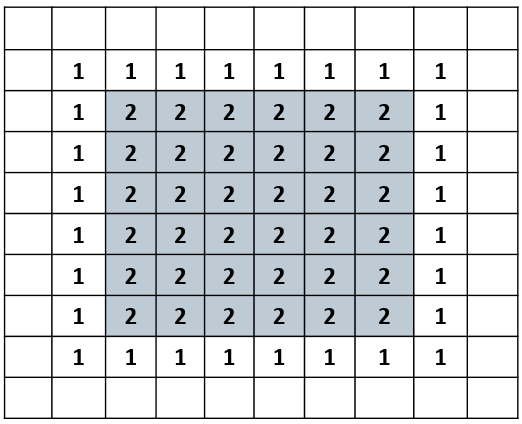}
	\caption{Computation (6x6)}
	\end{subfigure}
	\quad

	\begin{subfigure}[b]{0.3\textwidth}
	\includegraphics[width=1.0\textwidth]{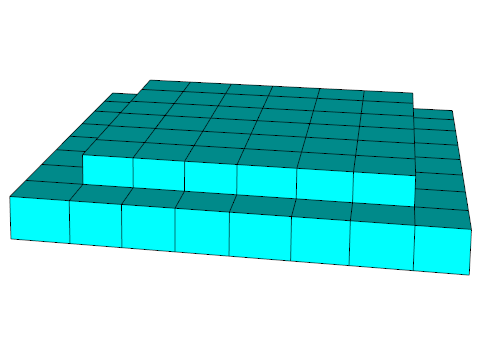}
	\caption{Level 1 (6x6)}
	\end{subfigure}
	\begin{subfigure}[b]{0.3\textwidth}
	\includegraphics[width=1.0\textwidth]{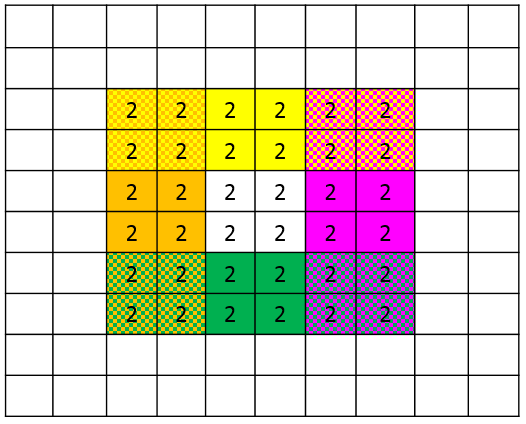}
	\caption{Populating the sides}
	\end{subfigure}
	\begin{subfigure}[b]{0.3\textwidth}
	\includegraphics[width=1.0\textwidth]{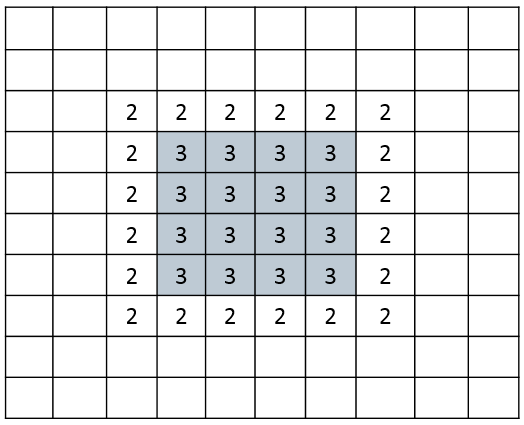}
	\caption{Computation (4x4)}
	\end{subfigure}
	\quad
	
	\begin{subfigure}[b]{0.3\textwidth}
	\includegraphics[width=1.0\textwidth]{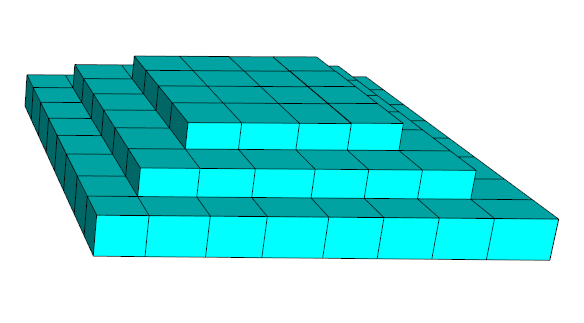}
	\caption{Level 2 (4x4)}
	\end{subfigure}
	\begin{subfigure}[b]{0.3\textwidth}
	\includegraphics[width=1.0\textwidth]{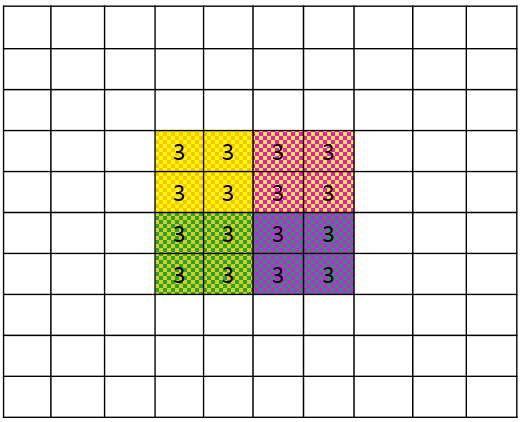}
	\caption{Populating the sides}
	\end{subfigure}
	\begin{subfigure}[b]{0.3\textwidth}
	\includegraphics[width=1.0\textwidth]{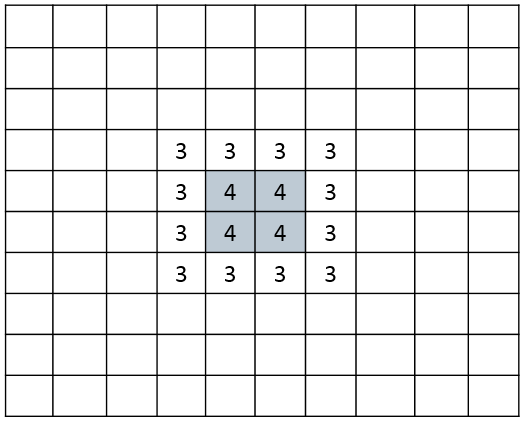}
	\caption{Computation (2x2)}
	\end{subfigure}
	\quad
	
	\begin{subfigure}[b]{0.45\textwidth}
	\includegraphics[width=1.0\textwidth]{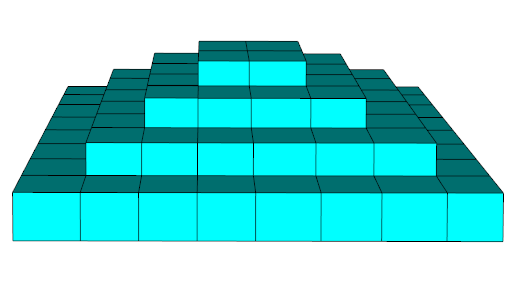}
	\caption{Level 3 (2x2)}
	\end{subfigure}
	\begin{subfigure}[b]{0.45\textwidth}
	\includegraphics[width=1.0\textwidth]{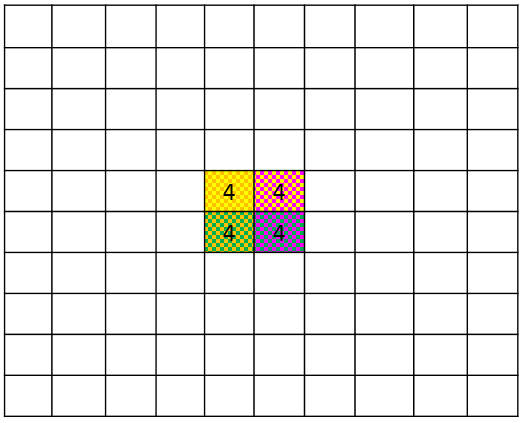}
	\caption{Populating the sides}
	\end{subfigure}

\caption{Illustrating the building process of the Upward Pyramid }
\label{upward}
\end{figure}

``Building'' the upward pyramid means computing the values at all the space-time grids
in the pyramid.  After the pyramid is built, the values at the space-time grids on the
four triangular sides, or upward pyramid panels, form the output, which feed into the other components of Swept 2D.  Building the pyramid
requires the values at the square base (time step $0$) as inputs, then applying the stencil
operation on the space-time grid points at time steps $1,0+2, \ldots, 0+n/2-1$.
.  As an example, if the square base is 8 by 8 and each grid point has an value
$u_{i,j}^{0}=1$; our stencil operation is incrementing by 1, i.e.,
$u_{i,j}^{k+1}=u_{i,j}^{k}+1$.  For illustration, we color code the 4 outputs,
the north, south, west, and east triangular sides, with yellow, green, orange,
and pink, respectively.\label{31}
Figure \ref{upward} illustrates the steps of building the upward pyramid.

The algorithm for building the upward pyramid is described in Algorithm \ref{upalgorithm}.  The algorithm has two $(n+2)$ by $(n+2)$ internal arrays, $\mathcal{U}$ and $\mathcal{D}$.\newline

\begin{algorithm}[H]
 \SetKwInOut{Input}{Input}
 \SetKwInOut{Output}{Output}
 $(\mathcal{N,S,W,E})$ = function \underline{UpwardPyramid} $(\mathbf{St},\mathcal{B})$\\
 \Input{$\mathbf{St}$: a list of stencil operations\\
        $\mathcal{B}$: an $n$ by $n$ array}
 \Output{$\mathcal{N,S,W,E}$: 4 arrays representing triangular sides}
 $\mathcal{D}_{1:n,1:n} \gets \mathcal{B}$\;
 $\mathcal{N} \gets \varnothing, \mathcal{S} \gets \varnothing, \mathcal{W} \gets \varnothing, \mathcal{E} \gets \varnothing$\;
 \For{$k=0,\ldots,\frac{n}{2}-1$}{
  $\mathcal{N}^k \gets \mathcal{D}_{k+1:n-k,k+1:k+2}$ \quad
  $\mathcal{N} \gets \mathcal{N} \cup \mathcal{N}^k $\\
  $\mathcal{S}^k \gets \mathcal{D}_{k+1:n-k,n-k-1:n-k}$ \quad
  $\mathcal{S} \gets \mathcal{S} \cup \mathcal{S}^k$\\
  $\mathcal{W}^k \gets \mathcal{D}_{k+1:k+2,k+1:n-k}$ \quad
  $\mathcal{W} \gets \mathcal{W} \cup \mathcal{W}^k$ \\
  $\mathcal{E}^k \gets \mathcal{D}_{n-k-1:n-k,k+1:n-k}$ \quad
  $\mathcal{E} \gets \mathcal{E} \cup \mathcal{E}^k$\\
 \For{$i=2+k,\ldots,n-k-1$}{
 \For{$j=2+k,\ldots,n-k-1$}{
 $\mathcal{U}_{i,j}
 = \mathbf{St}_k(\{\mathcal{D}_{i',j'},|i'-i|\le 1,|j'-j|\le 1\})$
 }}
 $\mathcal{U}\leftrightarrow\mathcal{D}$\\\
 }
 \caption{Building The Swept 2D Up-Pyramid}
 \label{upalgorithm}
\end{algorithm}
\qquad

Let us walk thought an example that has a square base of size 8 by 8.  The Upward pyramid algorithm will proceed as follows.  First, populate the 4  triangular sides, upward pyramid panels, with 2 layers from the four sides of the base.  As shown in Figure \ref{upward}(b), each triangular side now has 16 values.  Figure \ref{upward}(c) illustrates the subsequent stencil operation on level 0.

Proceeding to the next level, we further populate the 4 panels, each with 12 values from the new level, as shown in Figure \ref{upward}(d,e).  Now 12 new values will be added to each panel.  Figure \ref{upward}(f) illustrates the subsequent stencil operation on level 1.

Figure \ref{upward}(g-i) illustrates the process on the next level, in which 8 more values are added to each panel.  Figure \ref{upward}(j-k) illustrates the process on the final level, in which the same 4 values are added to each panel, and no further computation is possible. 

\subsection{Longitudinal and Latitudinal Bridges}

These two components differ only in their orientation.  We hereby refer to both as bridges.  The Swept 2D bridge is a three-dimensional, discrete space-time, structure.  Similar to the Swept 2D upward pyramid component, the two spatial dimensions are discretized with a grid indexed by $(i,j)$ and the time dimension is discretized with time steps indexed by $k$.  

\begin{figure}
	\centering
	\begin{subfigure}[b]{0.3\textwidth}
	\includegraphics[width=1.0\textwidth]{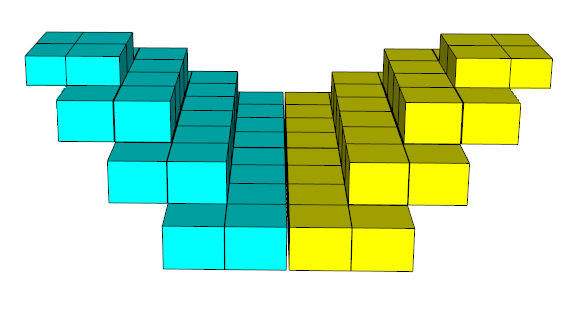}
	\caption{Linking two panels}
	\end{subfigure}
	\begin{subfigure}[b]{0.3\textwidth}
	\includegraphics[width=1.0\textwidth]{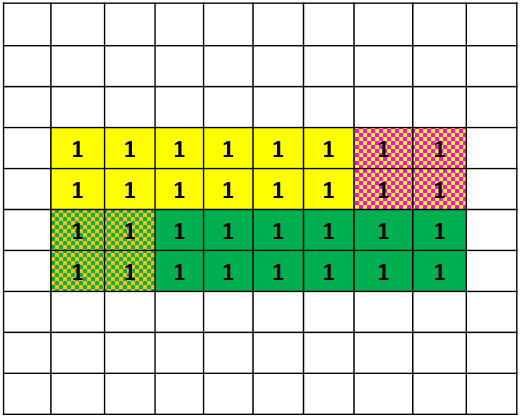}
	\caption{Populating panels}
	\end{subfigure}
	\begin{subfigure}[b]{0.3\textwidth}
	\includegraphics[width=1.0\textwidth]{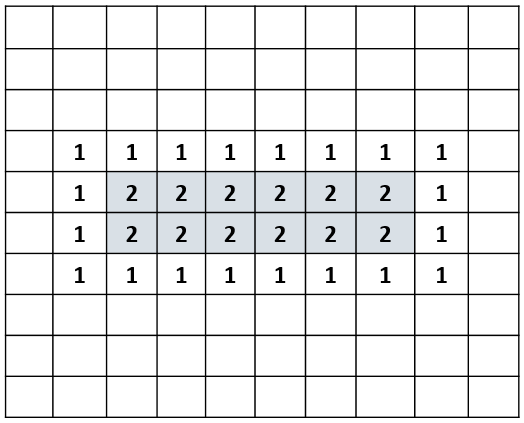}
	\caption{Computation}
	\end{subfigure}
	\quad

	\begin{subfigure}[b]{0.3\textwidth}
	\includegraphics[width=1.0\textwidth]{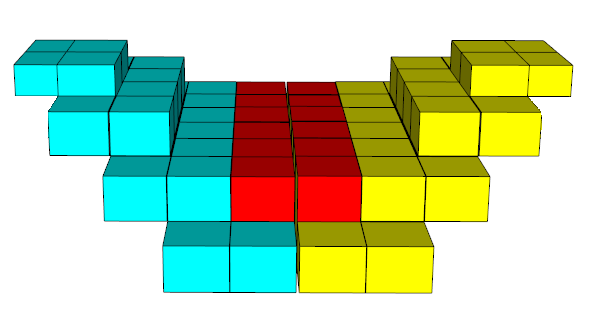}
	\caption{Moving Up in time}
	\end{subfigure}
	\begin{subfigure}[b]{0.3\textwidth}
	\includegraphics[width=1.0\textwidth]{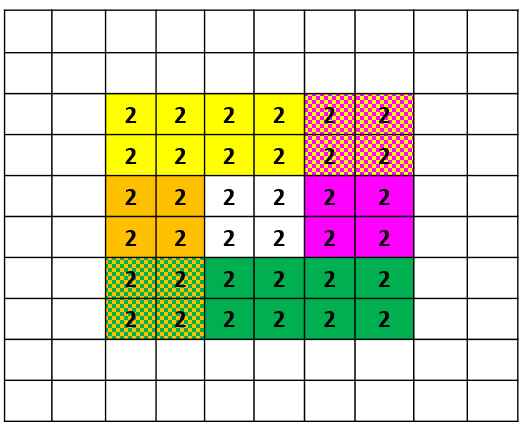}
	\caption{Populating panels}
	\end{subfigure}
	\begin{subfigure}[b]{0.3\textwidth}
	\includegraphics[width=1.0\textwidth]{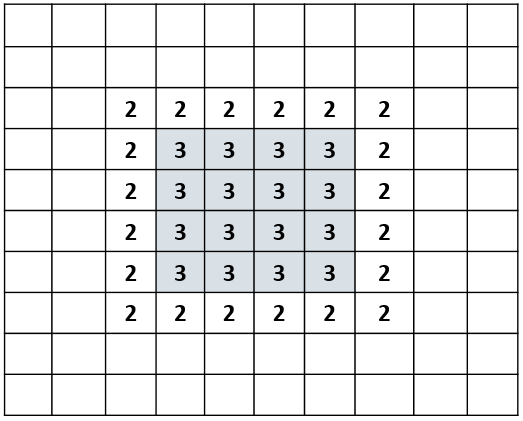}
	\caption{Computation}
	\end{subfigure}
	\quad
	
	\begin{subfigure}[b]{0.3\textwidth}
	\includegraphics[width=1.0\textwidth]{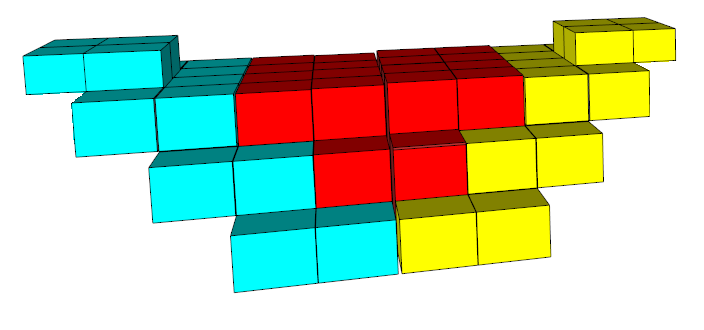}
	\caption{Moving Up in time}
	\end{subfigure}
	\begin{subfigure}[b]{0.3\textwidth}
	\includegraphics[width=1.0\textwidth]{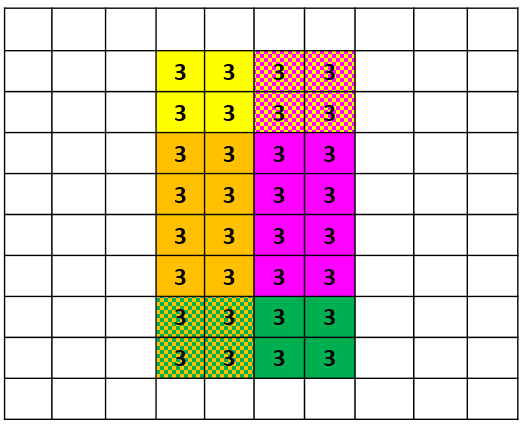}
	\caption{Populating panels}
	\end{subfigure}
	\begin{subfigure}[b]{0.3\textwidth}
	\includegraphics[width=1.0\textwidth]{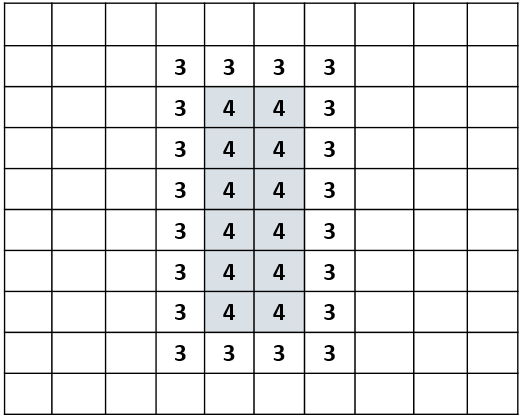}
	\caption{Computation}
	\end{subfigure}
	\quad
	
	\begin{subfigure}[b]{0.45\textwidth}
	\includegraphics[width=1.0\textwidth]{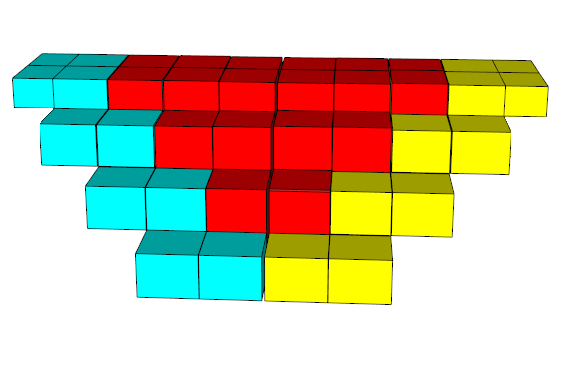}
	\caption{Reaching final level}
	\end{subfigure}
	\begin{subfigure}[b]{0.45\textwidth}
	\includegraphics[width=1.0\textwidth]{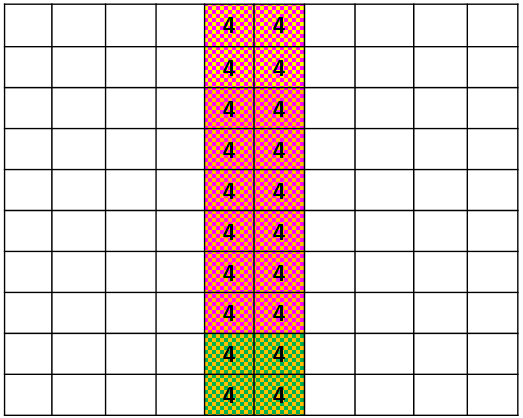}
	\caption{Populating panels}
	\end{subfigure}
	\caption{Illustrating the building process of the Longitudinal bridge}
	\label{vbridge}
\end{figure}

\begin{figure}
	\centering
	\begin{subfigure}[b]{0.3\textwidth}
	\includegraphics[width=1.0\textwidth]{./b1}
	\caption{Linking two panels}
	\end{subfigure}
	\begin{subfigure}[b]{0.3\textwidth}
	\includegraphics[width=1.0\textwidth]{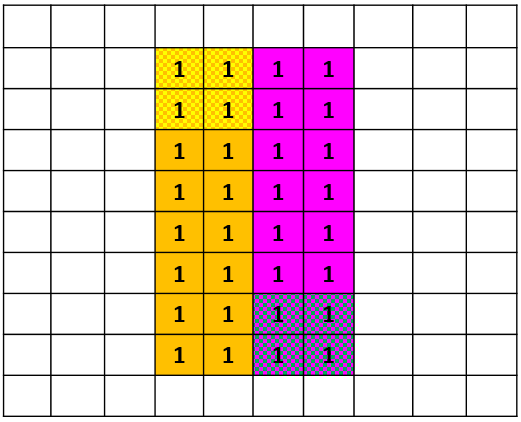}
	\caption{Populating panels}
	\end{subfigure}
	\begin{subfigure}[b]{0.3\textwidth}
	\includegraphics[width=1.0\textwidth]{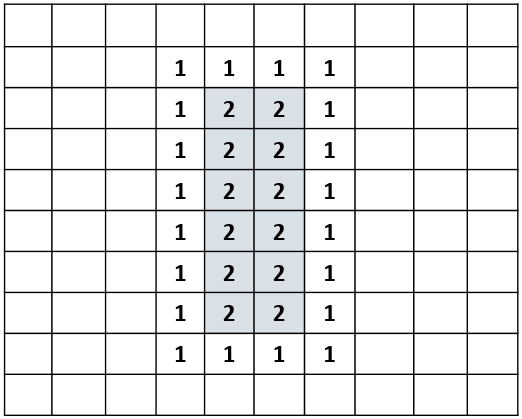}
	\caption{Computation}
	\end{subfigure}
	\quad

	\begin{subfigure}[b]{0.3\textwidth}
	\includegraphics[width=1.0\textwidth]{./b2}
	\caption{Moving Up in time}
	\end{subfigure}
	\begin{subfigure}[b]{0.3\textwidth}
	\includegraphics[width=1.0\textwidth]{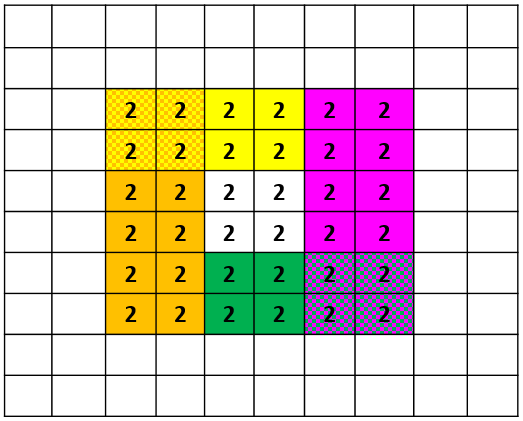}
	\caption{Populating panels}
	\end{subfigure}
	\begin{subfigure}[b]{0.3\textwidth}
	\includegraphics[width=1.0\textwidth]{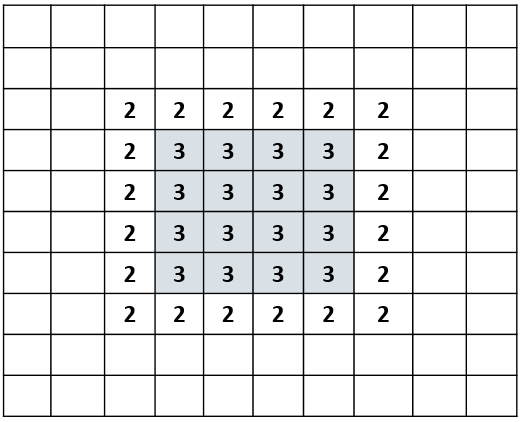}
	\caption{Computation}
	\end{subfigure}
	\quad
	
	\begin{subfigure}[b]{0.3\textwidth}
	\includegraphics[width=1.0\textwidth]{./b3}
	\caption{Moving Up in time}
	\end{subfigure}
	\begin{subfigure}[b]{0.3\textwidth}
	\includegraphics[width=1.0\textwidth]{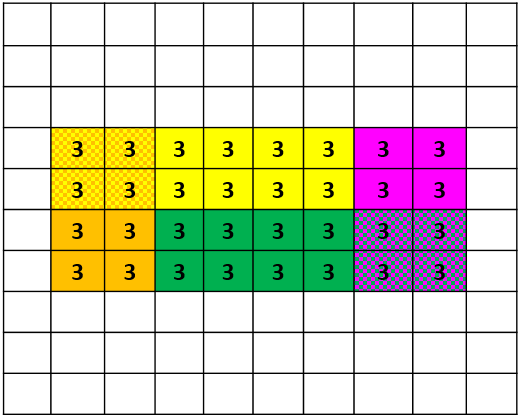}
	\caption{Populating panels}
	\end{subfigure}
	\begin{subfigure}[b]{0.3\textwidth}
	\includegraphics[width=1.0\textwidth]{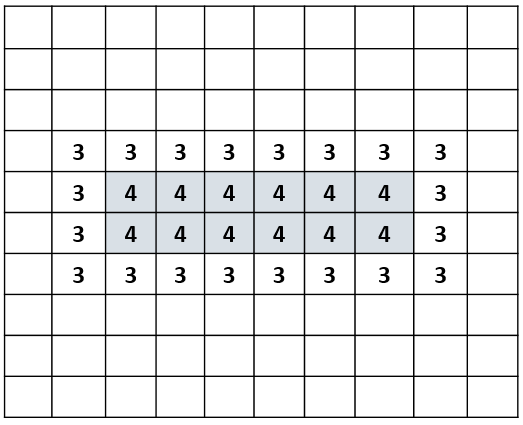}
	\caption{Computation}
	\end{subfigure}
	\quad
	
	\begin{subfigure}[b]{0.45\textwidth}
	\includegraphics[width=1.0\textwidth]{./b4}
	\caption{Reaching final level}
	\end{subfigure}
	\begin{subfigure}[b]{0.45\textwidth}
	\includegraphics[width=1.0\textwidth]{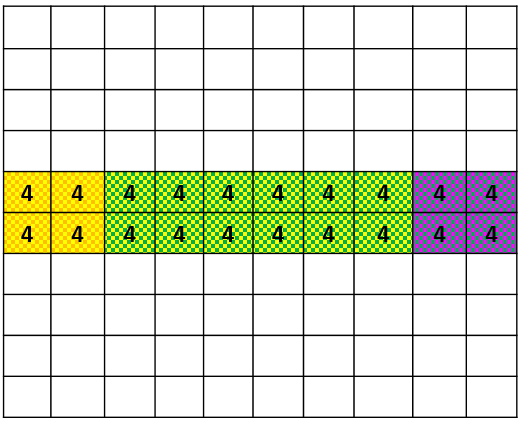}
	\caption{Populating panels}
	\end{subfigure}
	\caption{Illustrating the building process of the Latitudinal bridge}
	\label{hbridge}
\end{figure}

The Swept 2D bridges have the same height as the upward pyramid.  Thinking of the three-dimensional discrete space-time, the bridge fills a valley that resides between two adjacent upward pyramids.  ``Building'' the Swept 2D bridge means calculating all the space-time values in the valley, starting from 2 triangular sides, or panels, from two adjacent upward pyramids.  Figures \ref{vbridge} and \ref{hbridge} visualize how the bridge component is built.\newline

Depending on the orientation of the gap filled between the panels, we call the constructed bridge a Longitudinal or a Latitudinal bridge.  One may think of the difference between the Longitudinal and Latitudinal bridges as if the building process has loops with ``i'' and ``j'' indices that are basically flipped.\newline

The outputs of the bridge construction are two triangular sides.  They have the same shape and size of the upward pyramid sides, except that they are flipped in the time axis.  Their wider end is at the top and smaller end at the bottom.  A longitudinal bridge requires North and South upward pyramid sides and generates West and East bridge sides.  A latitudinal bridge requires West and East upward pyramid sides and produces North and South bridge sides.\newline

The algorithms for building the Longitudinal and Latitudinal bridges are described in Algorithms \ref{valgorithm} and \ref{halgorithm} respectively.  Just like the upward pyramid algorithm, both \ref{valgorithm} and \ref{halgorithm} algorithms have two $(n+2)$ by $(n+2)$ internal arrays, $\mathcal{U}$ and $\mathcal{D}$. Here $n$ is the side length of the square subdomain partition, the base of the upward pyramid.  If we denote the first time step in the Swept 2D bridge as $0$, the last time step will be $n/2-1$.

\begin{minipage}[t]{7cm}
\smaller

\begin{algorithm}[H]
 \SetKwInOut{Input}{Input}
 \SetKwInOut{Output}{Output}
 $(\mathcal{W,E})$ = function \underline{LongitudinalBridge} $(\mathbf{St},\mathcal{N,S})$\\
  \Input{$\mathbf{St}$: a list of stencil operations\\
         $\mathcal{N,S}$: 2 arrays representing triangular sides}
 \Output{$\mathcal{W,E}$: 2 arrays representing triangular sides}
 $\mathcal{W} \gets \varnothing, \mathcal{E} \gets \varnothing$\\
 \For{$k=0,\ldots,\frac{n}{2}-1$}{
 $\mathcal{D}_{k+1:n-k,\frac{n}{2}-k-1:\frac{n}{2}-k} \gets \mathcal{N}^k$ \\
 $\mathcal{D}_{k+1:n-k,\frac{n}{2}+k+1:\frac{n}{2}+k+2} \gets \mathcal{S}^k$ \\

 $\mathcal{W}^k \gets \mathcal{D}_{k+1:k+2,\frac{n}{2}-k+1:\frac{n}{2}+k+2}$ \quad
 $\mathcal{W} \gets \mathcal{W} \cup \mathcal{W}^k$ \\
 
 $\mathcal{E}^k \gets \mathcal{D}_{n-k-1:n-k,\frac{n}{2}-k-1:\frac{n}{2}+k}$ \quad
 $\mathcal{E} \gets \mathcal{E} \cup \mathcal{E}^k$\\
 
   \For{$i=2+k,\ldots,n-k-1$}{
 \For{$j=\frac{n}{2}-k,\ldots,\frac{n}{2}-k+1$}{
 $\mathcal{U}_{i,j}
 = \mathbf{St}_k(\{\mathcal{D}_{i',j'},|i'-i|\le 1,|j'-j|\le 1\})$
 }}
 $\mathcal{U}\leftrightarrow\mathcal{D}$\\
 }
 \caption{Building The Swept 2D Longitudinal Bridge}
 \label{valgorithm}
\end{algorithm}
\end{minipage} %
\begin{minipage}[t]{7cm}
\smaller

\begin{algorithm}[H]
 \SetKwInOut{Input}{Input}
 \SetKwInOut{Output}{Output}
 $(\mathcal{N,S})$ = function \underline{LatitudinalBridge} $(\mathbf{St},\mathcal{W,E})$\\
  \Input{$\mathbf{St}$: a list of stencil operations\\
         $\mathcal{W,E}$: 2 arrays representing triangular sides}
 \Output{$\mathcal{N,S}$: 2 arrays representing triangular sides}
 $\mathcal{N} \gets \varnothing, \mathcal{S} \gets \varnothing$\\
 \For{$k=0,\ldots,\frac{n}{2}-1$}{

 $\mathcal{D}_{\frac{n}{2}-k-1:\frac{n}{2}-k,k+1:n-k} \gets \mathcal{W}^k$ \\
 $\mathcal{D}_{\frac{n}{2}+k+1:\frac{n}{2}+k+2,k+1:n-k} \gets \mathcal{E}^k$ \\

 $\mathcal{N}^k \gets \mathcal{D}_{\frac{n}{2}-k+1:\frac{n}{2}+k+2,k+1:k+2}$ \quad
 $\mathcal{N} \gets \mathcal{N} \cup \mathcal{N}^k$ \\
 
 $\mathcal{S}^k \gets \mathcal{D}_{\frac{n}{2}-k-1:\frac{n}{2}+k,n-k-1:n-k}$ \quad
 $\mathcal{S} \gets \mathcal{S} \cup \mathcal{S}^k$\\

 \For{$i=\frac{n}{2}-k,\ldots,\frac{n}{2}-k+1$}{ 
   \For{$j=2+k,\ldots,n-k-1$}{
 $\mathcal{U}_{i,j}
 = \mathbf{St}_k(\{\mathcal{D}_{i',j'},|i'-i|\le 1,|j'-j|\le 1\})$
 }}
 $\mathcal{U}\leftrightarrow\mathcal{D}$\\
 }
 \caption{Building The Swept 2D Latitudinal Bridge}
 \label{halgorithm}
\end{algorithm}
\end{minipage}

\vspace{10pt}

Figure \ref{vbridge} shows the building of a Longitudinal bridge for a square domain partition of size 8 x 8. Figure \ref{vbridge}(a) shows the inputs, two upward pyramid triangular sides, or panels.  Here we will use the same color codes given in Section \ref{31} for the different upward pyramid sides in our example.

In the first step, we start by linking the first level of North and South upward pyramids triangular sides as shown in Figure \ref{vbridge}(a).  After that we populate the bottom part of newly generated West and East bridge triangular sides as illustrated in Figure \ref{vbridge}(b) .  Now stencil operation can be done to the inner blocks as shown in Figure \ref{vbridge}(c).

On the next level, as shown in Figure \ref{vbridge}(d), we place the values contained in the second level of our input panels, the North and South Upward Pyramid panels, next to those we just computed in the previous step.  We then populate the second level of the West and East bridge panels, our output, and then perform the second computation.  Figures \ref{vbridge}(d-f) clarify this process.  \ref{vbridge}(g-i) illustrate a similar process done on the next level, after which we reach the final level where no further computation is possible. We just populate the last level of the East and West bridge triangular sides, as shown in Figure \ref{vbridge}(j,k).

As mentioned earlier, the Latitudinal bridge can be build in a similar way.  So, we just summaries the process of building the Latitudinal bridge in Figure \ref{hbridge}.

\subsection{Downward Pyramid}

\begin{figure}
	\centering
	\begin{subfigure}[b]{0.45\textwidth}
	\includegraphics[width=1.0\textwidth]{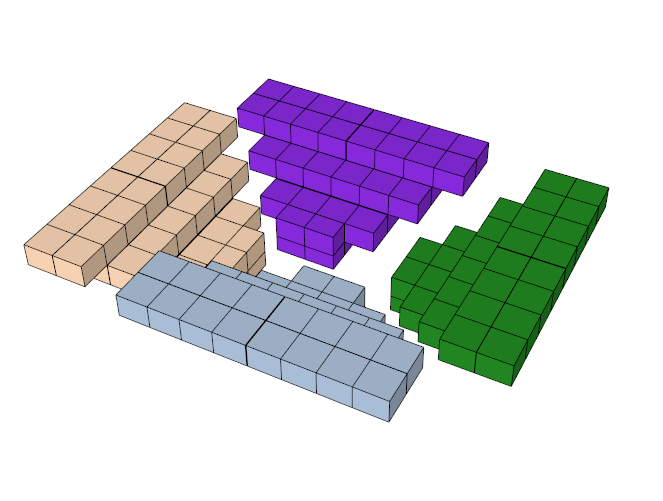}
	\caption{Staring with 4 bridges triangular sides}
	\end{subfigure}
	\begin{subfigure}[b]{0.45\textwidth}
	\includegraphics[width=1.0\textwidth]{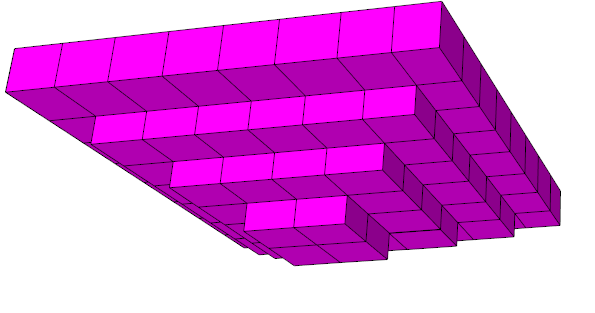}
	\caption{Computation is done to fill the gap}
	\end{subfigure}
	\caption{Abstracted view for building the Swept 2D Downward Pyramid. The computations in right figure fits into the gap between the triangular sides in the left figure.}
	\label{downabstract}
\end{figure}

This is the last component we need to complete the entire Swept Rule in 2D.  This component is also three dimensional in space-time.  Just like the Swept 2D upward pyramid and bridge components, the two spatial dimensions of the downward pyramid are discretized with a grid indexed by $(i,j)$ and the time dimension is discretized with time steps indexed by $k$.  Again, if we denote the first time step in the Swept 2D downward pyramid as $0$, the last time step will be $n/2$; where $n$ is the side length of the square subdomain partition. Note that the Swept 2D downward pyramids are one level higher in time than the upward pyramids and bridges components of Swept 2D.  ``Building'' the downward pyramid means calculating all the values in space-time of this three-dimensional component.  We construct a downward pyramid by filling a gap between four triangular sides of two Longitudinal and two Latitudinal bridges.  Figure \ref{downabstract}(a) shows the four triangular sides and the result obtained by filling the hollow area between them is shown in Figure \ref{downabstract}(b). The end result in \ref{downabstract}(b) is what we call the downward pyramid, which is a square pyramid that is similar in shape as the upward pyramid.  It is not difficult to see that the gap between the 4 triangular sides in Figure \ref{downabstract}(a) looks like a square pyramid that points downwards.  The algorithm for building the downward pyramid is described in Algorithm \ref{downalgorithm}.  The algorithm has two $(n+2)$ by $(n+2)$ internal arrays, $\mathcal{U}$ and $\mathcal{D}$.\newline

\begin{algorithm}[H]
 \SetKwInOut{Input}{Input}
 \SetKwInOut{Output}{Output}
 $(\mathcal{B})$ = function \underline{DownwardPyramid} $(\mathbf{St},\mathcal{N,S,W,E})$\\
  \Input{$\mathbf{St}$: a list of stencil operations\\
         $\mathcal{N,S,W,E}$: 4 arrays representing triangular sides}
 \Output{$\mathcal{B}$: an $n$ by $n$ array}

 \For{$k=0,\ldots,\frac{n}{2}$}{
 
 $\mathcal{D}_{\frac{n}{2}-k-1:\frac{n}{2}-k,\frac{n}{2}-k-1:\frac{n}{2}+k} \gets \mathcal{W}^k$ \\
 $\mathcal{D}_{\frac{n}{2}+k+1:\frac{n}{2}+k+2,\frac{n}{2}-k+1:\frac{n}{2}+k+2} \gets \mathcal{E}^k$ \\
 $\mathcal{D}_{\frac{n}{2}-k+1:\frac{n}{2}+k+2,\frac{n}{2}-k-1:\frac{n}{2}-k} \gets \mathcal{N}^k$ \\
 $\mathcal{D}_{\frac{n}{2}-k-1:\frac{n}{2}+k,\frac{n}{2}+k+1:\frac{n}{2}+k+2} \gets \mathcal{S}^k$ \\

 \For{$i=\frac{n}{2}-k,\ldots,\frac{n}{2}+k+1$}{ 
   \For{$j=\frac{n}{2}-k,\ldots,\frac{n}{2}+k+1$}{
 $\mathcal{U}_{i,j}
 = \mathbf{St}_k(\{\mathcal{D}_{i',j'},|i'-i|\le 1,|j'-j|\le 1\})$
 }}
 $\mathcal{U}\leftrightarrow\mathcal{D}$\\
 }
$\mathcal{B} \gets \mathcal{D}_{1:n,1:n}$\\
\qquad
\caption{Building The Swept 2D Downward Pyramid}
\label{downalgorithm}
\end{algorithm}
\quad

\begin{figure}
	\centering
	\begin{subfigure}[b]{0.4\textwidth}
	\includegraphics[width=1.0\textwidth]{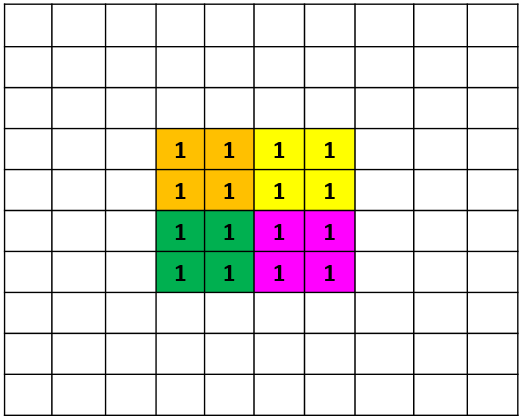}
	\caption{Linking the first level of 4 panels}
	\end{subfigure}
	\quad
	\begin{subfigure}[b]{0.4\textwidth}
	\includegraphics[width=1.0\textwidth]{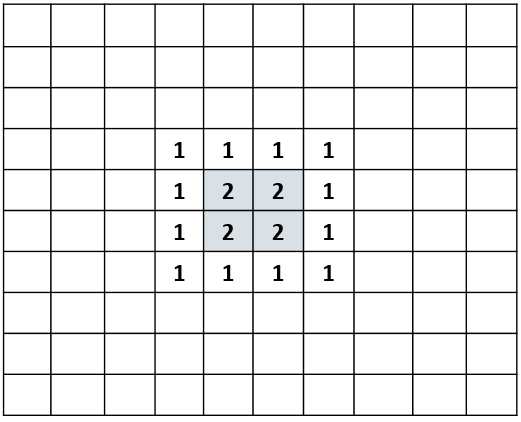}
	\caption{Computation}
	\end{subfigure}

	\begin{subfigure}[b]{0.4\textwidth}
	\includegraphics[width=1.0\textwidth]{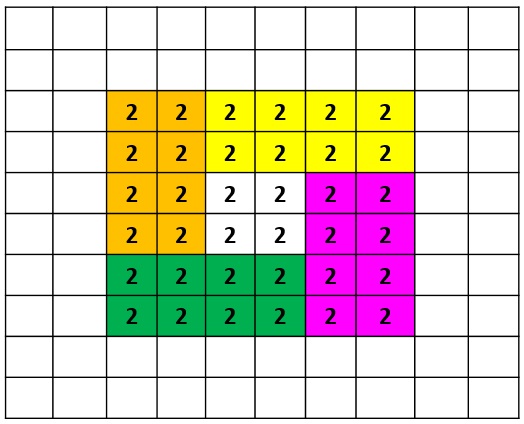}
	\caption{Placing the next level of 4 panels}
	\end{subfigure}
	\quad
	\begin{subfigure}[b]{0.4\textwidth}
	\includegraphics[width=1.0\textwidth]{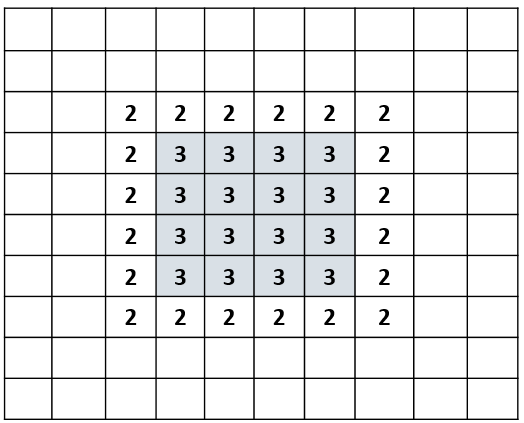}
	\caption{Computation}
	\end{subfigure}
	
	\begin{subfigure}[b]{0.4\textwidth}
	\includegraphics[width=1.0\textwidth]{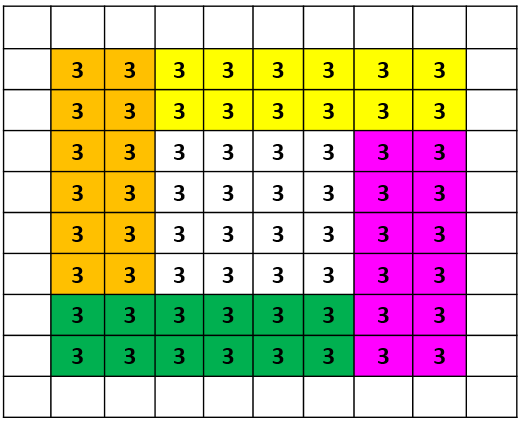}
	\caption{Placing the next level of 4 panels}
	\end{subfigure}
	\quad
	\begin{subfigure}[b]{0.4\textwidth}
	\includegraphics[width=1.0\textwidth]{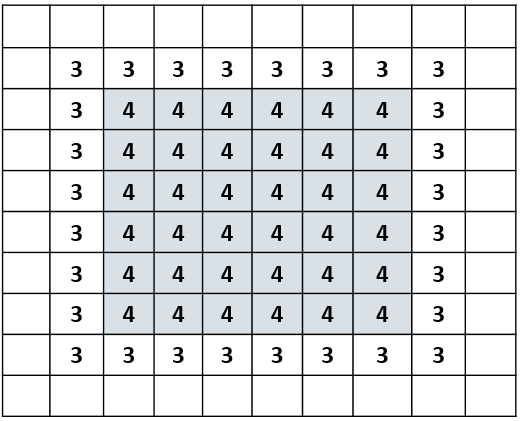}
	\caption{Computation}
	\end{subfigure}
	
	\begin{subfigure}[b]{0.4\textwidth}
	\includegraphics[width=1.0\textwidth]{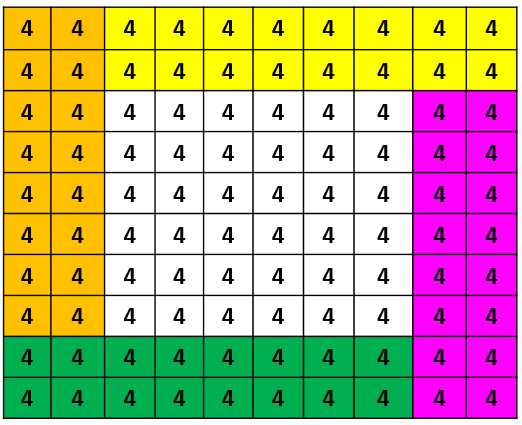}
	\caption{Placing the last level of 4 panels}
	\end{subfigure}
	\quad
	\begin{subfigure}[b]{0.4\textwidth}
	\includegraphics[width=1.0\textwidth]{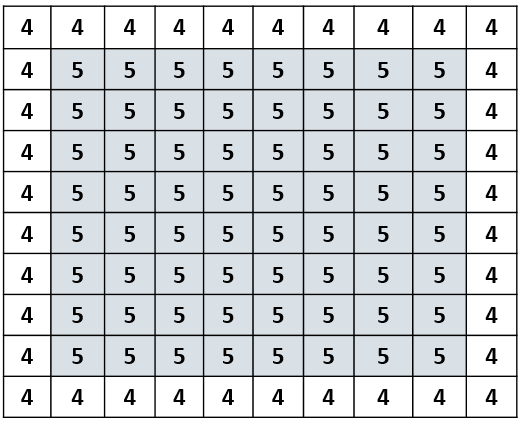}
	\caption{Computation}
	\end{subfigure}
    \caption{Illustrating the building process of the Downward Pyramid component of Swept 2D}
    \label{down}
\end{figure}

As an example, we build a downward pyramid starting from 4 bridge triangular sides, North, South, West, East, and our output will be the top level of a Downward Pyramid, which represents a solution of the PDE at the same time level.\newline

Starting with 4 panels, we simply place their first level values appropriately to form a 4 by 4 square as shown in Figure \ref{down}(a), where computation for the inner 2 by 2 square is possible.  Figure \ref{down}(b)clarifies the process.  Proceeding to the next level, we place the next level of our four inputs next to the results from the previous level to form a 6 by 6 square. That allows us to produce the inter 4 by 4 square of the next level, as clarified by Figures \ref{down}(c,d).  We proceed through the next levels until the computation is done to the entire 8 by 8 square subdomain as illustrated in Figures \ref{down}(e-h).

\section{Connecting the Swept Rule 2D components}

Notice how each block at any level of our 4 Swept 2D components has the complete set of immediate neighbors from the level below.  Let us relate that to the numerical stencil-based 2D simulation world.  Assume that each block is a grid point in finite difference, a controlled volume in finite volume, or an element in finite element discretization.  The blocks below each block represent the complete 9-point 2D stencil in a numerical discretization.\newline

\begin{figure}
\centering
\includegraphics[width=0.7\linewidth]{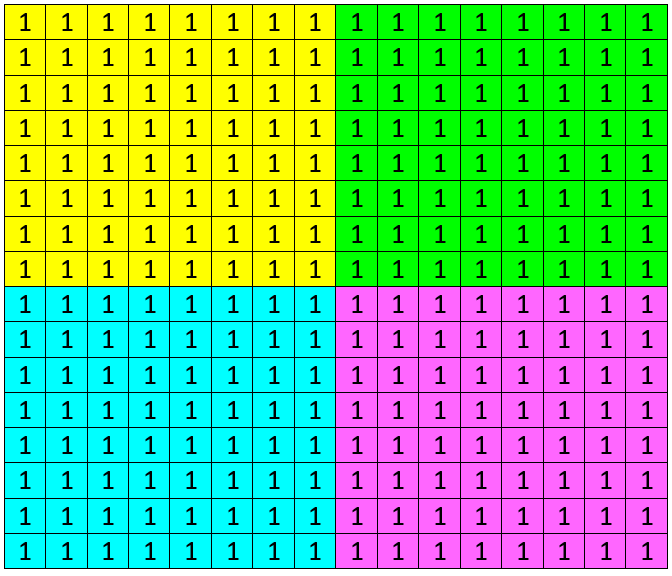}
\caption{Illustration of a square domain partitioned into 4 sub-domains.  Initial time level is 1}
\label{initial}
\end{figure}

After the visual clarification of our 4 Swept 2D components, it is now time to show how these components work with each other to build the Swept Rule in 2D.  Assume that we are starting with a square computational domain, with periodic boundary conditions, that can be decomposed into 4 subdomains.  Each of these subdomains is assigned to a different processor.  For simplicity in illustrating how the components are built, we will show a top view of our square domain with numbers representing the level of each grid point, i.e. timestep, of each block.  Figure \ref{initial} illustrates our square domain decomposed into 4 smaller square subdomains.\newline

\begin{figure}
\centering
\begin{subfigure}[b]{0.4\textwidth}
	\includegraphics[width=1.0\textwidth]{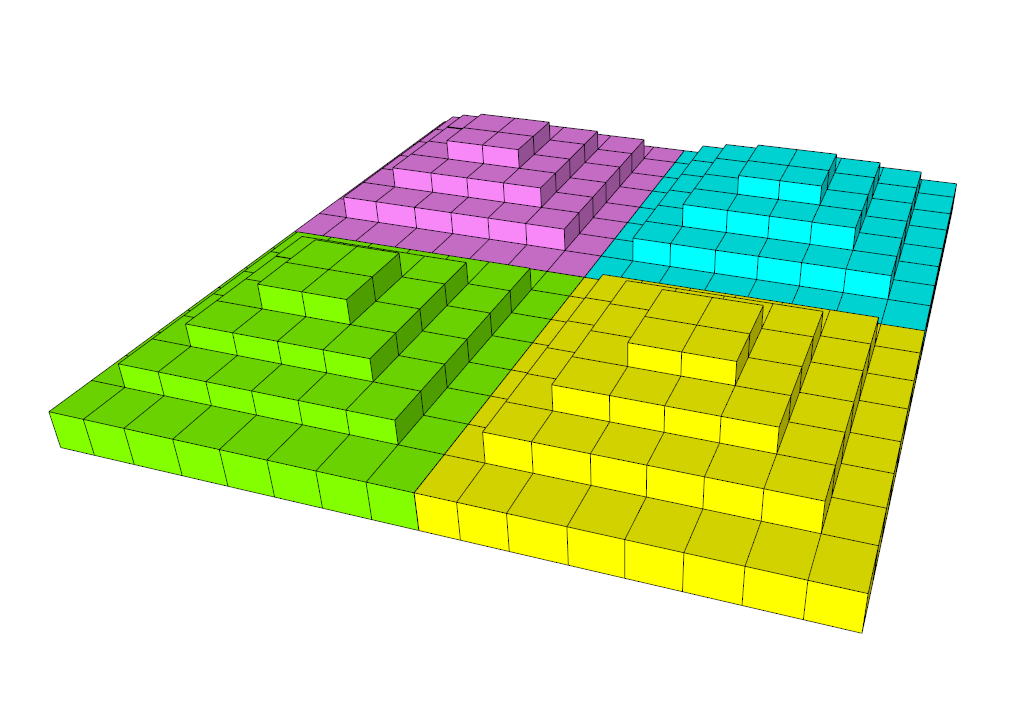}
	\caption{The square domain state after Swept 2D stage 1 is complete}
	\end{subfigure}
	\quad
	\begin{subfigure}[b]{0.4\textwidth}
	\includegraphics[width=1.0\textwidth]{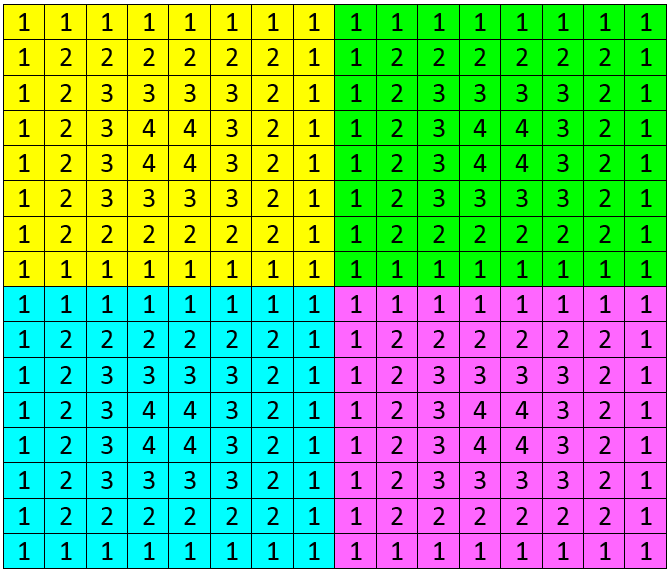}
	\caption{A top-view of the square domain}
	\end{subfigure}

\caption{The square domain state after Swept 2D stage 1 is complete}
\label{firstdone}
\end{figure}

As each processor starts to explicitly solve the PDE, following the domain of influence the domains of dependency to progress without the need to communicate with neighbouring processes, each processor will be building an upward pyramid.  Figure \ref{firstdone}(a) shows the end result of the first stage of the Swept Rule in 2D.  A 2D top view of the end result of stage 1 is shown in Figure \ref{firstdone}(b) with numbers representing the time level of each block.  Notice how the blocks are at different time levels.  The tip of the Upward Pyramids at each subdomain is at level 4.\newline

\begin{figure}
	\centering
	\begin{subfigure}[b]{0.45\textwidth}
	\includegraphics[width=1.0\textwidth]{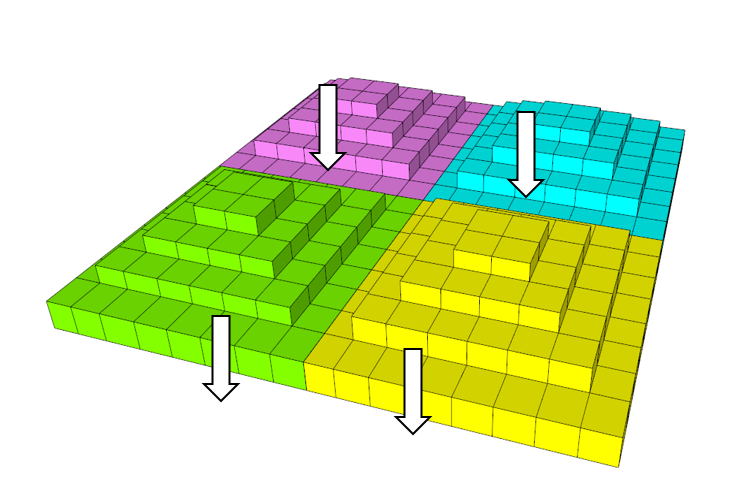}
	\caption{A logical illustration of the locations of the Longitudinal bridges}
	\end{subfigure}
	\quad
	\begin{subfigure}[b]{0.45\textwidth}
	\includegraphics[width=1.0\textwidth]{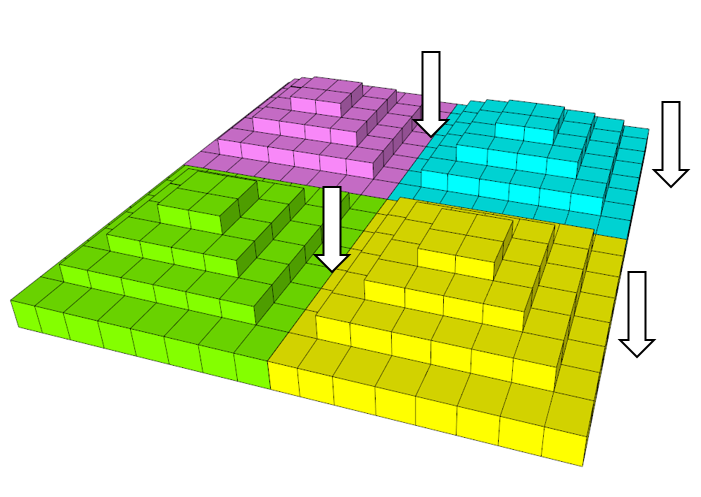}
	\caption{A logical illustration of the locations of the Latitudinal bridges}
	\end{subfigure}
		
    \caption{Locations of the Longitudinal and Latitudinal bridges of Swept 2D}
    \label{virtualbridges}
\end{figure}

The next stage in Swept 2D will be to build the Longitudinal and Latitudinal bridges.  Those will fill the Longitudinal and Latitudinal aligned gaps between two pyramids.  Figures \ref{virtualbridges}(a,b) show the locations of the Longitudinal and Latitudinal bridges to be built.  Notice that the arrows pointing outside the domain represent the bridges that are to be build due to the periodic boundary conditions.  Basically, there will be a set of Longitudinal and Latitudinal bridges between the boundary pyramids.  We use this technique to avoid working in a fraction, i.e. quarter or half, Swept 2D component.\newline

\begin{figure}
\centering
	\begin{subfigure}[b]{0.45\textwidth}
	\includegraphics[width=0.7\linewidth]{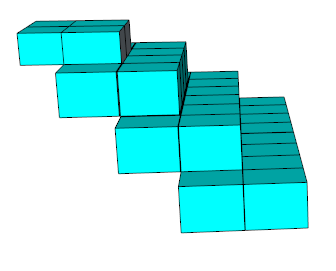}
	\caption{The shape of an Upward Pyramid triangular side (panel)}
	\end{subfigure}
	\begin{subfigure}[b]{0.45\textwidth}
	\includegraphics[width=0.7\linewidth]{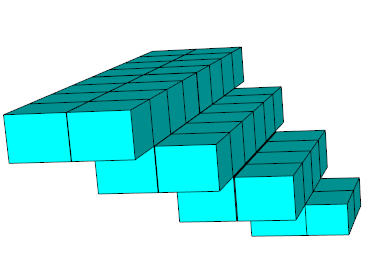}
	\caption{The shape of a Bridge triangular side (panel)}
	\end{subfigure}
\quad
\caption{The shapes of upward pyramid and bridges triangular sides}
\label{panels}
\end{figure}

In order to proceed from stage 1 to stage 2, that is building the bridges, the first communication between the processors needs to take place.  Each process needs to send data to two of its neighbors and receive data from its other two neighbors.  To be more specific, each process will exchange 2 triangular sides, panels, of the upward pyramid that it built in stage 1.  Figure \ref{panels}(a) clarifies what we mean by an upward pyramid triangular side or panel.\newline

If we agree that each upward pyramid has 4 triangular side (panels) and those panels can be named as North, South, West and East, then we notice that in order to build a Longitudinal bridge, each process needs a set of North and South panels.  On the other hand, to build a Latitudinal bridge, each process needs a set of West and East panels.  For this reason, each process will communicate its East and South upward pyramid panels to its East and South neighbors respectively.  Each process will also receive 2 upward pyramid panels from its North and West neighbors.  With these panels, each process can build a pair of bridges; one Longitudinal and one Latitudinal bridge.\newline

Notice that building the bridges at the boundary will cause a shift of the entire computational domain due to the periodic boundary condition assumed at the beginning.  Again, this is to avid splitting a single Swept 2D component between processors.  This way, each process will work on a complete component throughout the Swept 2D process .\newline

The last stage will be to fill the remaining gaps between the bridges.  This is the step where the Downward Pyramids are built.  From our previous explanation of building the downward pyramid, we see that 4 triangular sides (panels) are needed.  The panels this time come from the bridges that we built in stage 2.  Notice that each Longitudinal bridge generates West and East panels.  Also, each Latitudinal bridge generates North and South panels.  The panels generated by the bridges are similar to those previously generated by the Upward pyramids, except that they are flipped in the time axis.  Figure \ref{panels}(b) illustrates bridge panels.

So, in order to proceed with building the Downward Pyramids, a second communication needs to take place.  The communication this time is also for exchanging 2 bridge panels.  Again, each process will send its East and South bridges panels to its East and South neighbors respectively.  Each process will also receive 2 panels from its North and West neighbors.  After the panels of the bridges are properly exchanged between the processes, each can build a downward pyramid.\newline

\begin{figure}
	\centering
	\begin{subfigure}[b]{0.45\textwidth}
	\includegraphics[width=1.0\textwidth]{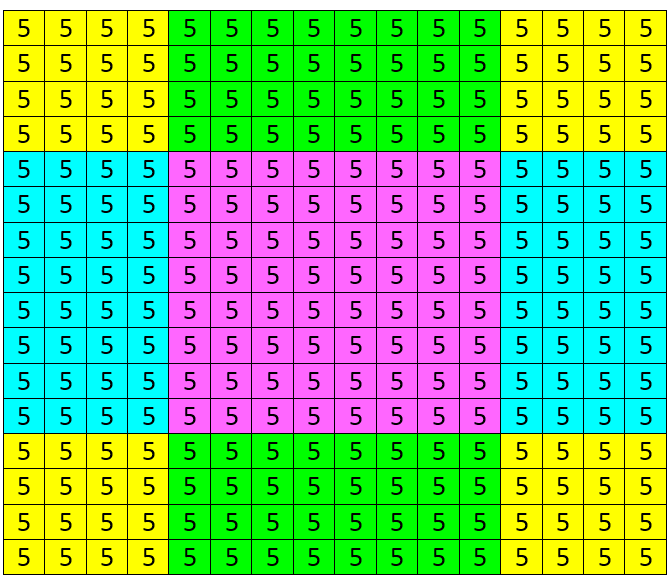}
	\caption{The state of square domain after a half Swept 2D cycle}
	\end{subfigure}
	\quad
	\begin{subfigure}[b]{0.45\textwidth}
	\includegraphics[width=1.0\textwidth]{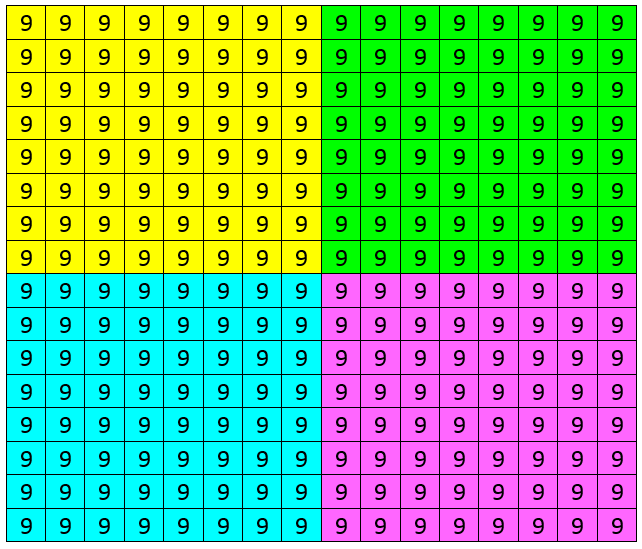}
	\caption{The state of square domain after a complete Swept 2D cycle}
	\end{subfigure}
		
    \caption{Showing the state of square domain after a half and a complete Swept 2D cycle}
    \label{halffull}
\end{figure}

At the end of this stage, the entire computational domain is at a consistent state.  Meaning that all blocks, grid points or cells, are at the same time level.  However, the domain arrangement has changed as a result of the shift that took place due to the periodic boundary condition.  We call what we did so far a half Swept 2D cycle.  Figure \ref{halffull}(a) clarifies what we mean by a computational domain shift due to periodic boundary condition.  The view in Figure \ref{halffull}(a) is logical; the actual work was done to complete Swept 2D components.  Each color in the figure represents work done by a single process\newline

The other half of the Swept 2D cycle can be performed in the same way, except that the data exchange between the processes will happen between the other two neighbors.  For example, instead of sending the North and West triangular sides (panels) of the Upward Pyramid and bridges to North and West neighbors, we send the South and East triangular sides (panels) of upward pyramids and bridges to the South and East neighbors.  At the end of the second half of the Swept 2D cycle, the domain will have its original arrangement.  Figure \ref{halffull}(b) shows the computational domain after a complete Swept 2D cycle.\newline

So, a complete Swept 2D cycle, requires a total of 4 communications to take place.  2 communications happen at each half of the Swept 2D cycle.  Notice how starting with an initial condition as time level 1, and working with square subdomains of 8 by 8 allows a complete Swept 2D cycle to promote the computational domain 8 time levels involving only 4 communications between the processors.\newline

We showed earlier how the Swept 2D components are simple and easy to build.  Now we show the algorithm that connects, or glues, the Swept 2D components together making them construct the Swept Rule in 2D.  Going through Algorithm \ref{swept2d} clarifies many points that may not have been clarified.  One such point, is how we mange the communication of the panels of Upward Pyramids and Bridges.  For simplicity, we abbreviated the names of the panels that get exchanged between the processors by their direction and whether they are upward (belonging to an upward pyramid) or downward (belonging to a bridge).  For example, $\mathcal{E_U}$, stands for an East-Up panel and $\mathcal{N_D}$, stands for a North-Down panel.  Moreover, we used the ``$\circleddotright$'' symbol to denote sending to a process and the ``$\circleddotleft$'' to denote receiving from a process.  We also abbreviated the neighboring processes with ``$\mathbb{P}$'' symbol combined with the direction of that process.  For example, $\mathbb{P}_N$ means the North neighboring process.  So, in algorithm \ref{swept2d}, a line that reads $\mathcal{N_U} \circleddotright  \mathbb{P}_N[\mathcal{N_U}]$, simply means sending the North triangular side of the Upward Pyramid to the North process.\newline

\begin{algorithm}[H]
 \SetKwInOut{Input}{Input}
 \SetKwInOut{Output}{Output}

$(\mathcal{B})$ = function \underline{ Swept 2D} $(\mathbf{St},\mathcal{B},\mathcal{C})$\\
 \Input{$\mathbf{St}$: a list of stencil operations\\
        $\mathcal{B}$: an $n$ by $n$ array representing a square subdomain\\
        $\mathcal{C}$: The number of Swept 2D cycles to perform}
 \Output{$\mathcal{B}$: The square subdomain partition after applying $\mathcal{C}*n$ subtimesteps}
 
 \For{$c=1,\ldots,\mathcal{C}$}{
  
  $\mathcal{N_U,S_U,W_U,E_U}\gets $ \underline{UpwardPyramid} $(\mathbf{St},\mathcal{B})$\\
  $\mathcal{N_U} \circleddotright  \mathbb{P}_N[\mathcal{N_U}] ;\qquad   \mathcal{W_U}  \circleddotright \mathbb{P}_W[\mathcal{W_U}] $\\ 
  $\mathcal{N_U} \circleddotleft   \mathbb{P}_S[\mathcal{N_U}] ;\qquad   \mathcal{W_U}  \circleddotleft  \mathbb{P}_E[\mathcal{W_U}] $\\
  $\mathcal{W_D,E_D} \gets $  \underline{LongitudinalBridge}   $(\mathbf{St},\mathcal{S_U,N_U})$\\
  $\mathcal{N_D,S_D} \gets $  \underline{LatitudinalBridge} $(\mathbf{St},\mathcal{E_U,W_U})$\\
  $\mathcal{N_D} \circleddotright  \mathbb{P}_N[\mathcal{N_D}] ;\qquad   \mathcal{W_D}  \circleddotright \mathbb{P}_W[\mathcal{W_D}] $\\ 
  $\mathcal{N_D} \circleddotleft   \mathbb{P}_S[\mathcal{N_D}] ;\qquad   \mathcal{W_D}  \circleddotleft  \mathbb{P}_E[\mathcal{W_D}] $\\
  $\mathcal{B} \gets $ \underline{DownwardPyramid} $(\mathbf{St},\mathcal{S_D,N_D,E_D,W_D})$\\
  \vspace{1em}
  
  $\mathcal{N_U,S_U,W_U,E_U}\gets $ \underline{UpwardPyramid} $(\mathbf{St},\mathcal{B})$\\
  $\mathcal{S_U} \circleddotright  \mathbb{P}_S[\mathcal{S_U}] ;\qquad   \mathcal{E_U}  \circleddotright \mathbb{P}_E[\mathcal{E_U}] $\\ 
  $\mathcal{S_U} \circleddotleft   \mathbb{P}_N[\mathcal{S_U}] ;\qquad   \mathcal{E_U}  \circleddotleft  \mathbb{P}_W[\mathcal{E_U}] $\\
  $\mathcal{W_D,E_D} \gets $  \underline{LongitudinalBridge}   $(\mathbf{St},\mathcal{S_U,N_U})$\\
  $\mathcal{N_D,S_D} \gets $  \underline{LatitudinalBridge} $(\mathbf{St},\mathcal{E_U,W_U})$\\  
  $\mathcal{S_D} \circleddotright  \mathbb{P}_S[\mathcal{S_D}] ;\qquad   \mathcal{E_D}  \circleddotright \mathbb{P}_E[\mathcal{E_D}] $\\ 
  $\mathcal{S_D} \circleddotleft   \mathbb{P}_N[\mathcal{S_D}] ;\qquad   \mathcal{E_D}  \circleddotleft  \mathbb{P}_W[\mathcal{E_D}] $\\
  $\mathcal{B} \gets $ \underline{DownwardPyramid} $(\mathbf{St},\mathcal{S_D,N_D,E_D,W_D})$\\
  
 }
 \caption{The Swept 2D}
 \label{swept2d}
\end{algorithm}

%% file: analysis.tex
\section{A simplified performance analysis of the swept rule}

To simplify our analysis, let us consider a half Swept 2D cycle.  We will also assume that if we partition the computational domain into small square subdomains, bandwidth issues are not encountered during communications.  Meaning that, Communication between computing nodes takes time $\tau$, regardless of how much data is communicated.\newline

Now, if we partition our computational domain into square sub-domains of side length $n$, then within a process, each Swept 2D half cycle will perform  $n/2$ substeps for $n^2$ grid points, which is the area of each subdomain.  We know that each Swept 2D half cycle requires 2 communications between the computation processes.  Let us assume that each physical compute node contains a single MPI process, then for a half Swept 2D cycle we will encounter 2 communication latencies between our compute nodes.\newline

Let us break the half Swept 2D cycle to its basic components and try to understand how the 2 communication latencies are distributed among these components.  In half a cycle, we build one Upward Pyramid, one Longitudinal Bridge, one and Latitudinal Bridge and finally a Downward Pyramid.\newline

Building each of the 4 basic components involves some overhead that is a function of $n$ (assuming linear relationship for small values of n).  Let us denote each of these overhead values as $\alpha_u(n)$ for Upward Pyramids, $\alpha_d(n)$ for Downward Pyramids, and 
$\alpha_b(n)$ for both Longitudinal and Latitudinal Bridges.  Also, assume that the communication latency between the physical compute nodes, regardless of how much data is exchanged, is $\tau$ and the time to perform each sub-timestep per grid point is $s$.\newline

Knowing that half Swept 2D cycle performs $n/2$ sub-timesteps for $n^2$ grid points, we can calculate the time needed to perform each sub-timestep as follows:

\begin{equation}
[(n^2 s) \frac{n}{2} + \alpha_u(n) + \alpha_d(n) + 2 \alpha_b(n) + 2 \tau]/\frac{n}{2}
\label{fullcost}
\end{equation}

If we ignore the overhead cost values and simplify equation \ref{fullcost}, we get:

\begin{equation}
n^2s + \frac{4\tau}{n}
\label{cost}
\end{equation}

Now we present a simplified theoretical performance model for the Swept Rule in 2D similar to that we presented for the Swept Rule in a single space dimension\cite{swept}.  To understand how fast the swept rule in 2D is, we need to know the typical values of $\tau$ and $s$.\newline

\begin{table}
\smaller
\begin{center}
\begin{tabular}{|l|l|}\hline

Interconnect                   & Typical latency ($\tau$) \\\hline
Amazon EC2 cloud	           & 150 $\mu s$               \\\hline
Typical Gigabit Ethernet	   & 50  $\mu s$               \\\hline
Fast 100-Gigabit Ethernet	   & 5   $\mu s$               \\\hline
Mellanox 56Gb/s FDR InfiniBand & .7  $\mu s$               \\\hline

\end{tabular}
\end{center}
\caption{The range of latency $\tau$ commonly encountered today.}
\end{table}

Table 1 attempts to cover the range of latency $\tau$ one may encounter today.  The latency can change over three orders of magnitude, from the fastest Infiniband to a cloud computing environment not designed for PDEs.
\newline
The range of $s$ is even wider; it can span over eight orders of magnitude. $s$ depends both on the computing power of each node and on the complexity of each sub-timestep. If one sub-timestep on one spatial point takes $f$ floating point operations (FLOP) to process, then processing an array of them with a node capable of $F$ floating point operations per second (FLOPS) takes $s = f/F$ seconds per step-point. Running a cheap discretization on a powerful computing node leads to small $f$ and large $F$, therefore a small $s$; running an expensive discretization on a light node leads to large $f$ and small $F$, therefore a large $s$.\newline

Table 2 attempts to estimate the range of $s$ by covering the typical $f$ for solving PDEs in a two spatial dimensions, as well as the highest and lowest $F$ on a modern computing node. Consider $f$ of the heat equation, discretized with finite difference. Each sub-timestep takes only few FLOPs. A sub-timestep of a nonlinear system of equations, discretized with high order finite element, may take thousands of FLOPs. The highest $F$ is achieved today by GPU nodes. The Summit supercomputer, expected to be delivered to Oak Ridge Leadership Computing Facility (OLCF) in 2017, uses a similar architecture to achieve 40 TeraFLOPS per node. Using older CPU architecture is slower. Particularly slow is equivalence of using a single thread per node, e.g., in flat-MPI-style parallel programming. On the first generation Intel Core i3/i5/i7 architecture, codenamed Nehalem, a single thread is capable of about 10 GFLOPS. These different node architectures, running different PDE discretization schemes, yield orders of magnitude different values of $s$.\newline

\begin{table}
\smaller
\smaller
\begin{center}
\begin{tabular}{|c|c|c|}\hline

Computing node -- FLOPS & FLOP per step-point & Computing time\cr & & per step-point ($s$) \\\hline
Single thread Intel Nehalem & & \cr 10 GFLOPS & 8000 (FE system) & 800 $ns$ \\\hline
Single thread Intel Nehalem & & \cr 10 GFLOPS & 400 (FV system)  & 40 $ns$  \\\hline
Single thread Intel Nehalem & & \cr 10 GFLOPS & 6 (FD scalar)    & 0.6 $ns$ \\\hline
Oak Ridge 2017 Summit node & & \cr  40 TFLOPS & 8000 (FE system) & 200 $ps$ \\\hline
Oak Ridge 2017 Summit node & & \cr  40 TFLOPS & 400 (FV system)  & 10 $ps$   \\\hline
Oak Ridge 2017 Summit node & & \cr  40 TFLOPS & 6 (FD scalar)    & 150 $fs$  \\\hline

\end{tabular}
\end{center}
\caption{The typical time it takes to compute one sub-timestep on a single spatial point for solving PDEs in 2 spatial dimensions}
\end{table}

\begin{figure}
\centering
\includegraphics[scale=0.5]{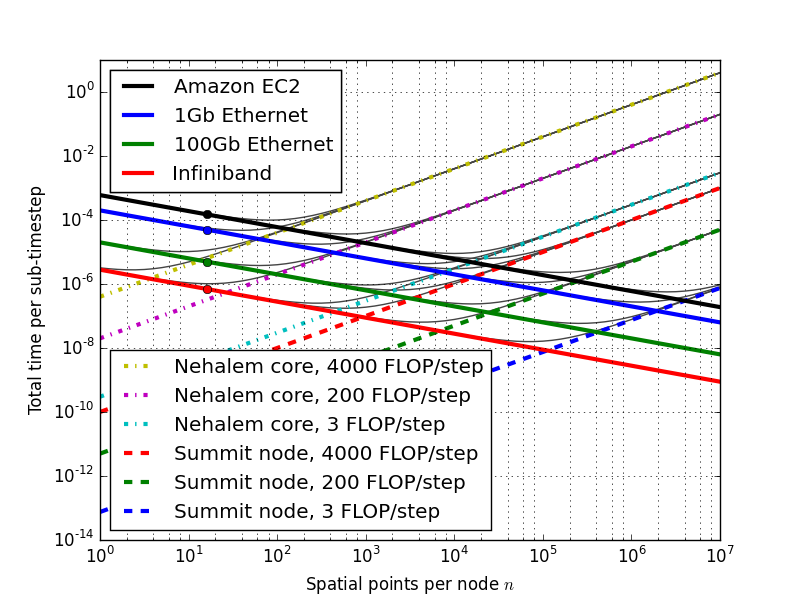}
\caption{Analyzing the Swept 2D performance using communication latency and CPU FLOP/step}
\label{trend}
\end{figure}

With these values of $\tau$ and $s$, the plot in Figure \ref{trend} shows, according to our sub-timespte cost formula and ignoring overhead cost, how fast the Swept 2D scheme runs as a function of $n$, the spatial points per node. The up-sloping, dashed and dash-dot lines represent the $n^2s$ term in our cost formula for different values of $s$, and the down-sloping, solid lines represent the $\frac{4\tau}{n}$ term for different values of $\tau$. For each combination of $s$ and $\tau$, the total time per sub-timestep as a function of $n$, plotted as thin black curves in Figure \ref{trend}, can be found by summing the corresponding up-sloping and down-sloping lines.\newline

The total time per sub-timestep can be minimized by choosing $n$. This minimizing $n$ can be found, for each $\tau$ and $s$.  The optimal value of $n$ represent the limit of scaling. Above this optimum, decreasing $n$ by scaling to more nodes would decrease the total time per sub-timestep, accelerating the simulation. But decreasing $n$ beyond the optimum by scaling to even more nodes would not accelerate, but slow down the simulation.\newline

At the optimal $n$, the swept rule 2D almost always breaks the latency barrier. A method that requires communication every sub-timestep takes at least $\tau$ per sub-timestep. This limit is the latency barrier.

%% file: interface.tex
\section{Interface and implementation of the swept scheme}

The interface we present for our implementation of the Swept Rule in 2D is very simple.  It is similar to that we presented earlier for the Swept Rule in single space-dimension\cite{swept}.  We basically present a data initialization function and a time-stepping function.

The input variables to the initialization function are the global ``i'' and ``j'' indicies for each spacial point and a spacial point structure representing a 2D stencil. The following pseudo code exemplifies the initialization function interface:

\begin{lstlisting}
Init(pointIndex i,pointIndex j,spacialPoint2D p)
{
    Based on the location of the point (i,j), set its initial data
	p.inputs[0] = variable 0 initial value
	p.inputs[1] = variable 1 initial value
	p.inputs[n] = variable n initial value
}
\end{lstlisting}

The second function in our interface is where the PDE solve is performed.  The input variables to the timestepping function are the index of which sub-timestep to be executed and a 2D 9-point stencil spacial point structure. The following psedocode shows the timestepping function of our interface.

\begin{lstlisting}
timestep(timestep i,spacialPoint2D p)
{
    Based on the value of i, perform the proper operation on p

}
\end{lstlisting}

We standardized our implementations of the Swept Rule on this interface so that it is easy to keep improving on the internal implementation of the Swept 2D components, while not modifying the numerical aspects of the implemented PDE solver.  It is know that during the development of numerical simulation applications, the concerns of numerical engineers are usually very different from those of software engineers or computer scientists.  a C++ implementation of the Swept Rule in 2D can be found at:\url{https://github.com/hubailmm/Wave2D} 

%% file: experiments.tex
\section{Swept rule in 2D -- solution of the 2D Wave Equation}

\begin{figure}
	\centering
	\begin{subfigure}[b]{0.3\textwidth}
	\includegraphics[width=1.0\textwidth]{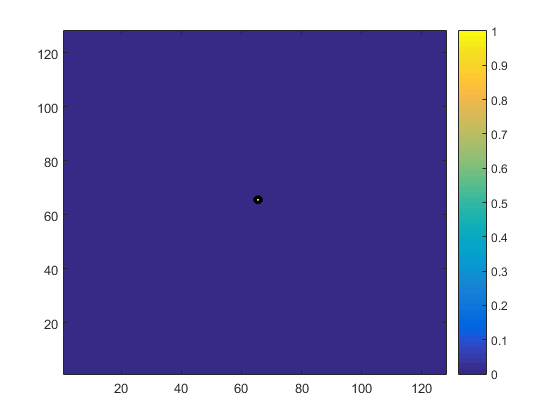}
	\caption{Intial condition}
	\end{subfigure}
	\begin{subfigure}[b]{0.3\textwidth}
	\includegraphics[width=1.0\textwidth]{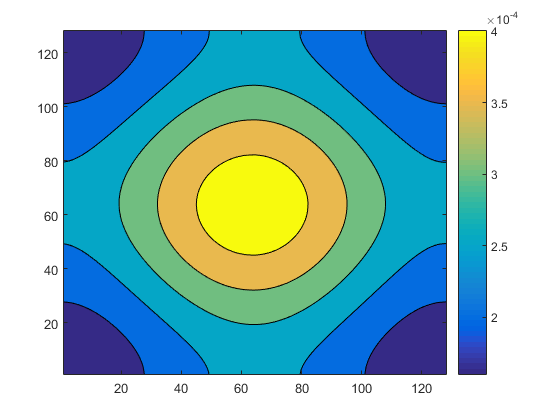}
	\caption{1280 timesteps}
	\end{subfigure}
	\begin{subfigure}[b]{0.3\textwidth}
	\includegraphics[width=\textwidth]{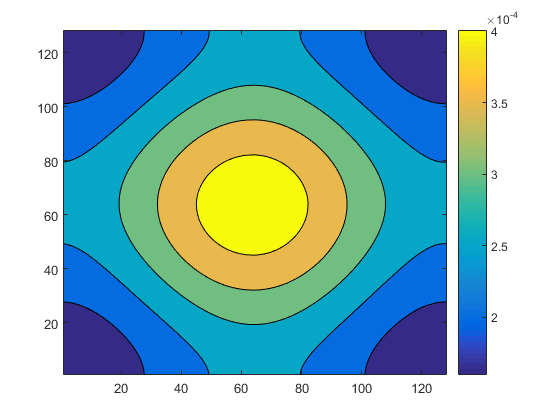}
	\caption{2560 timesteps}
	\end{subfigure}
	
	\quad
	
	\begin{subfigure}[b]{0.3\textwidth}
	\includegraphics[width=\textwidth]{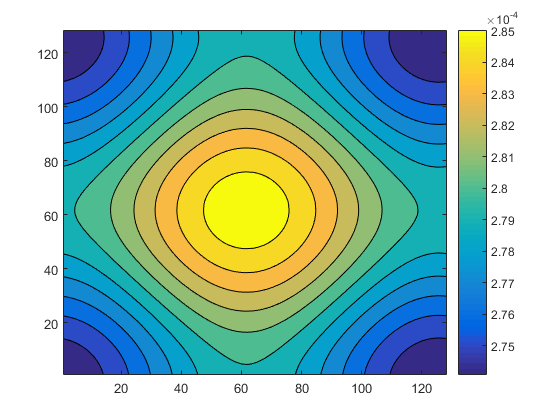}
	\caption{3840 timesteps}
	\end{subfigure}
	\begin{subfigure}[b]{0.3\textwidth}
	\includegraphics[width=\textwidth]{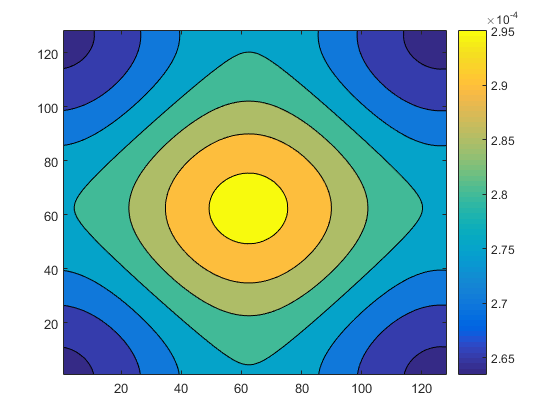}
	\caption{5120 timesteps}
	\end{subfigure}
	\begin{subfigure}[b]{0.3\textwidth}
	\includegraphics[width=\textwidth]{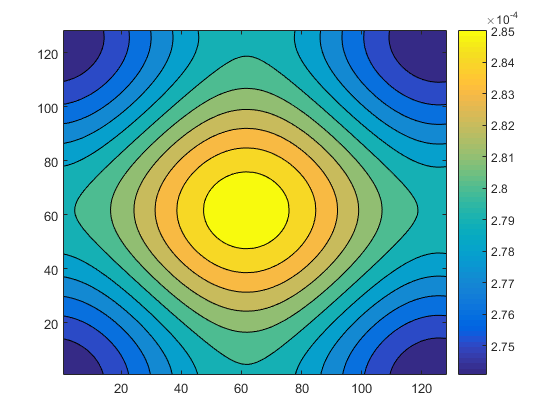}
	\caption{6400 timesteps}
	\end{subfigure}	
		
    \caption{Contour plots showing the propagation of the Wave.  This verifies our Swept 2D implementation}
    \label{wave}
\end{figure}

To start with, we tested our implementation of the Swept Rule in 2D in solving the two-dimensional wave equation.  Our PDE configuration was periodic boundary conditions, $CFL$ Number equals $0.3$ and an initial wave source located at the center of the domain.  The verification of our solution is shown in Figure \ref{wave}.  

An experiment was conducted in a small 9-node Amazon EC2 cluster.  StartCluster was used to form the 9-node cluster with an EC2 instance type of ``c3.large" and an Amazon Machine Image (AMI) of ``ami-6b211202"\cite{starcluster}.  A single MPI process was assigned to each node.  Figure \ref{waveperformance} plots in log scale the performance comparison between the straight ``Classical'' and Swept domain decomposition.  Notice how the Swept Rule in 2D gains a speed-up factor of more than 3 when we assign about 1000 grid points per MPI process.   Our implementation of the Swept Rule based solver of the 2D Wave equation can be found at:\url{https://github.com/hubailmm/Wave2D} 

\begin{figure}
\centering
\includegraphics[width=0.7\linewidth]{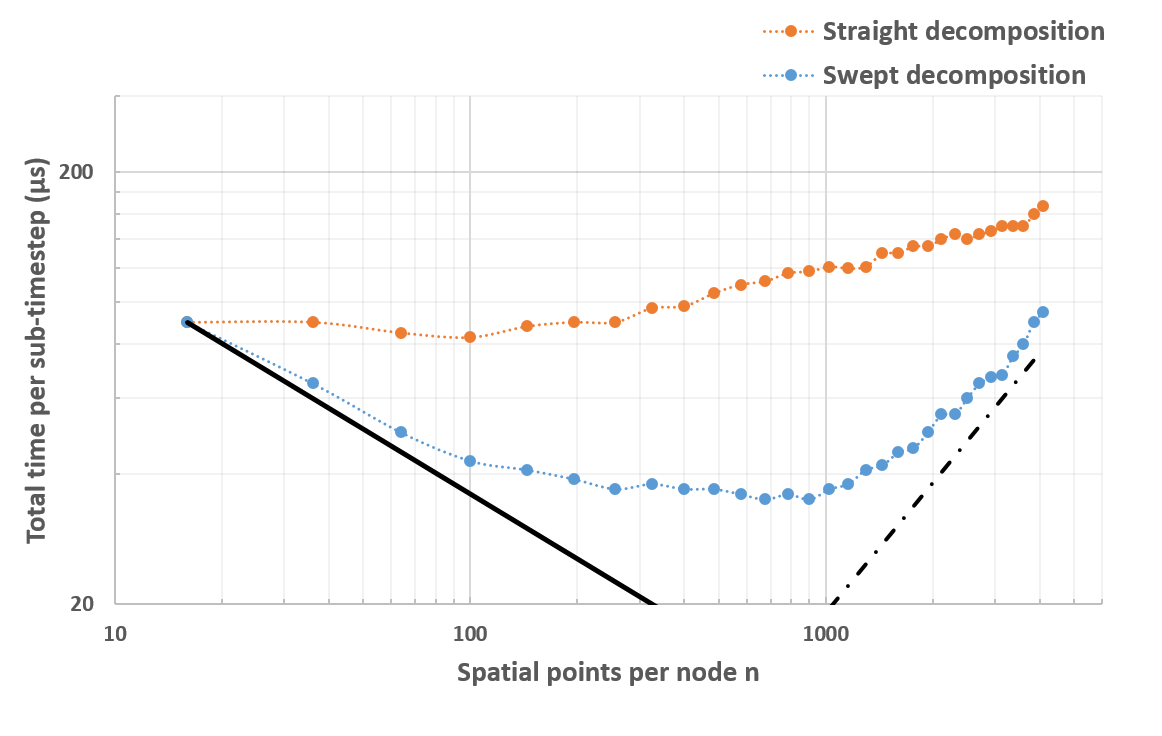}
\caption{The performance of Straight and Swept 2D domain partitioning schemes when solving the 2D wave equation.  The x-axis represents the number of spacial points assigned to each process.  The y-axis represents the time, in microseconds, needed to perform a timestep.  The black up-sloping, dashed line represents the $n^2s$ term in Equation \ref{cost}, and the black down-sloping, solid line represents the $ \frac{4\tau}{n}$ term in the Equation}
\label{waveperformance}
\end{figure}

\section{Swept rule in 2D -- solution of the 2D Euler Equation for Gas-Dynamics}

The Swept Rule in 2D was also tested on the Euler Equation for Gas-Dynamics with conservative, skew-symmetric, and compressible flow with periodic boundary conditions\cite{euler2d}.  

In the experiment, we used Finite Difference spacial discretization along with a 4-stage Runge-Kutta time integration scheme.  The experiment setup was a rectangular wind tunnel with $lx = 50 , ly = 25 , nx = 1024 , ny = 512 , dx = \frac{lx}{nx} , dy = \frac{ly}{ny} , dt = 1e-6$ and an obstacle, given by $e^{(x^2 + (y+.25*\frac{ly}{nx})^2)^8}$.  The simulation initial condition was $\rho = 1.084 , M = 0.2 , u = c * M, v = 0.0 , p = 101325$\newline

\begin{figure}
\centering
	\begin{subfigure}[b]{0.45\textwidth}
	\includegraphics[width=1.0\textwidth]{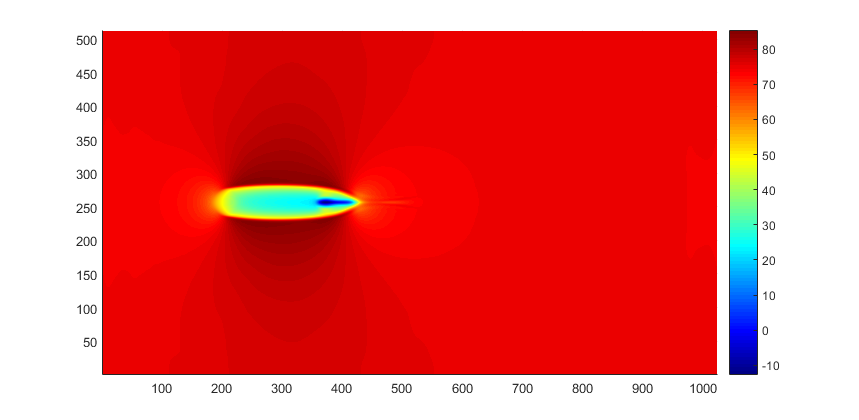}
	\caption{3200 timesptes}
	\end{subfigure}
	\begin{subfigure}[b]{0.45\textwidth}
	\includegraphics[width=\textwidth]{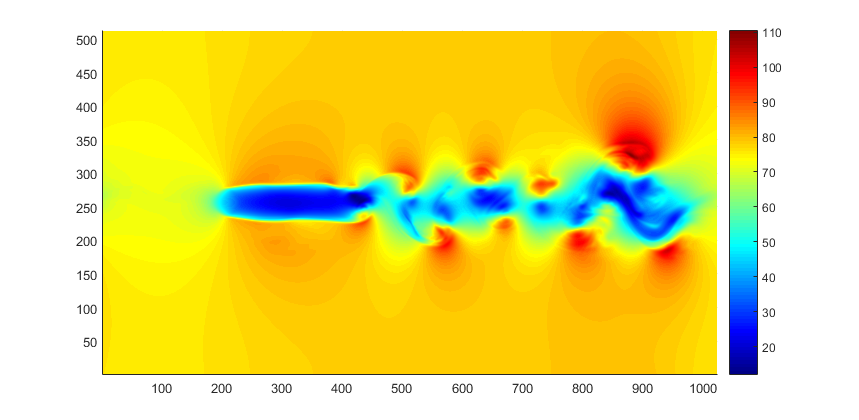}
	\caption{32000 timesteps}
	\end{subfigure}
\caption{Showing contour plot of the Horizontal Momentum (u) when solving the 2D Euler equation after 3200 and 32000 timesteps.  The X-axis and Y-axix represent the spacial location.}
\label{eulertest}
\end{figure}

The implementation was done in C++ and MPI.  Solution verification is shown in Figures \ref{eulertest}(a,b).  The Figures show the results obtained for the x momentum after 3200 and 32000 timesteps respectively.\newline

Our experiment was conducted in a cluster that was formed in Amazon's EC2 cloud services.  StartCluster was used to form the cluster with an EC2 instance type of ``c3.large" and an Amazon Machine Image (AMI) of ``ami-6b211202"\cite{starcluster}.  A single MPI process was assigned to each node.  Figure \ref{eulerperformance} compares the performance of the classic and Swept 2D partitioning schemes.  Notice how the Swept Rule in 2D achieves a speed-up factor of 4 when assigning about 400 grid points per node.  Our implementation of the Swept Rule based solver of the 2D Euler PDE can be found at:\url{https://github.com/hubailmm/Euler2D} 

\begin{figure}
\centering
\includegraphics[scale=0.4]{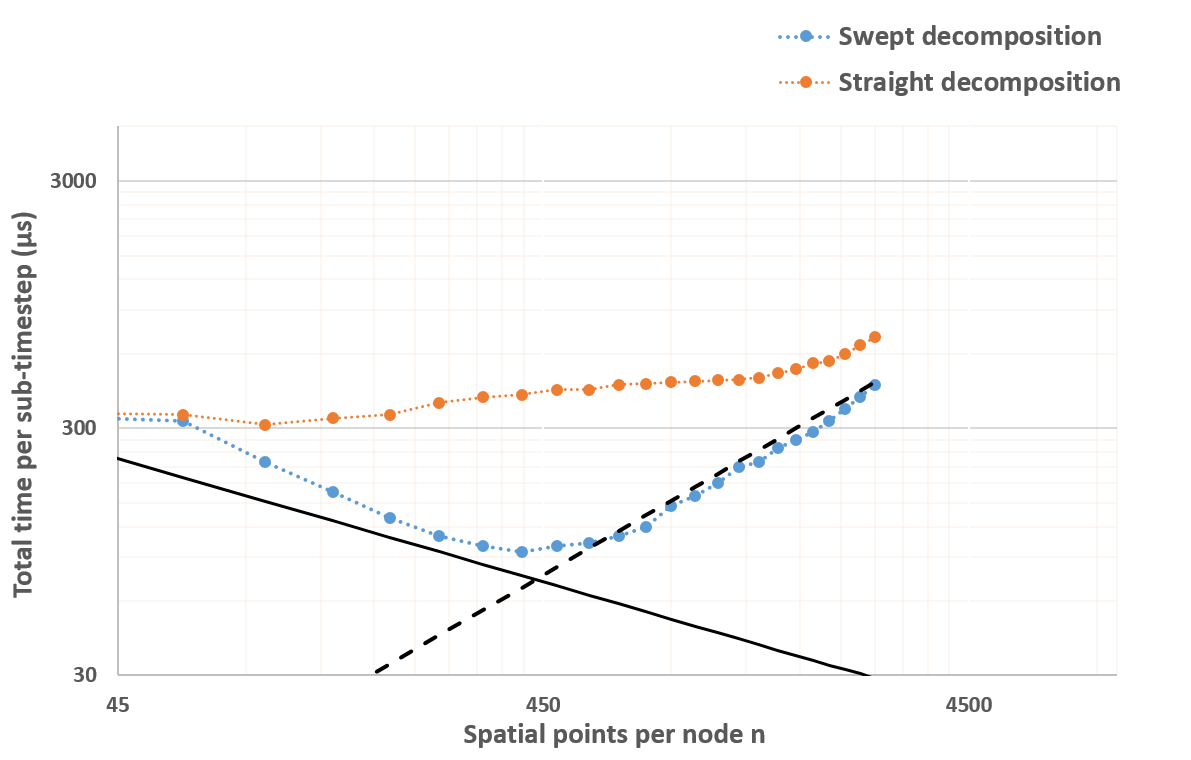}
\caption{The performance of Straight and Swept 2D domain partitioning schemes when solving the 2D Euler equation.  The x-axis represents the number of spacial points assigned to each process.  The y-axis represents the time, in microseconds, needed to perform a sub-timestep.  The black up-sloping, dashed line represents the $n^2s$ term in Equation \ref{cost}, and the black down-sloping, solid line represents the $ \frac{4\tau}{n}$ term in the Equation}
\label{eulerperformance}
\end{figure}

%% file: conclusion.tex
\section{Conclusion}

We presented the extension of the swept rule for decomposing space and time in solving
PDEs to two space-dimensional problems. The swept rule in 2D breaks the latency barrier, advancing each sub-timestep in a fraction of the time it takes for a message to travel from one computing
node to another. In our experiment with the 2D Wave and Euler equations on Amazon EC2, Swept 2D had a speed-up factor of more than 3 compared to the classical straight way of domain decomposition.

The Swept rule in 2D can be utilized by an interface that separates the concerns of the numerical scheme developer and Swept components developer.  The presented interface makes it very easy to develop different types of numerical application using the same internal implementation of the Swept 2D.  The authors' implementation of the swept decomposition scheme, along with testcases and illustrations in this article, are open source.\newline

\subsection*{Acknowledgment}
We acknowledge the Advanced Research Center at Saudi Aramco for sponsoring
graduate students to pursue research in the Computational Science and
Engineering area. We also acknowledge NASA NRA Award 15-TTT1-0057 under
Dr. Eric Nielsen and Dr. Mujeeb Malik, AFOSR Award F11B-T06-0007
under Dr. Fariba Fahroo and DOE award DE-FG02-14ER26173/DE-SC00011089 
under Dr. Sandy Landsberg.